\documentclass[11pt,a4paper]{article}
\usepackage[T1]{fontenc}
\usepackage[utf8]{inputenc}
\usepackage{lmodern}
\usepackage[margin=15mm]{geometry}
\usepackage{amsmath,amssymb,amsthm,mathtools}
\usepackage{microtype,enumitem,booktabs,array,tikz}
\usepackage[hidelinks]{hyperref}
\usepackage{cleveref}
\usepackage{cite}
\numberwithin{equation}{section}
\newtheorem{theorem}{Theorem}[section]
\newtheorem{proposition}[theorem]{Proposition}
\newtheorem{lemma}[theorem]{Lemma}

\theoremstyle{definition}

\theoremstyle{remark}

\newcommand{\R}{\mathbb R}
\newcommand{\C}{\mathbb C}

\newcommand{\dd}{\mathrm d}
\DeclareMathOperator{\rank}{rank}
\DeclareMathOperator{\diag}{diag}

\DeclareMathOperator{\corank}{corank}
\DeclareMathOperator{\adj}{adj}
\DeclareMathOperator{\Ran}{Ran}
\setlist{itemsep=2pt,topsep=4pt}
\hypersetup{pdftitle={Singular bifurcations in implicit discrete dynamical systems},pdfauthor={Jorge Buescu and Henrique M. Oliveira}}
\title{Singular bifurcations in implicit discrete dynamical systems}
\author{Jorge Buescu\thanks{\raggedright Department of Mathematical Sciences and CEMS.UL,
Faculty of Sciences, Universidade de Lisboa, Portugal.
\texttt{jsbuescu@ciencias.ulisboa.pt}.}
\and Henrique M. Oliveira\thanks{\raggedright Corresponding author.
Department of Mathematics and CAMGSD, Instituto Superior T\'ecnico,
Universidade de Lisboa, Av. Rovisco Pais, 1049--001 Lisbon, Portugal.
\texttt{henrique.m.oliveira@tecnico.ulisboa.pt}.}}
\date{}
\begin{document}
\maketitle
\begin{abstract}
A periodic orbit of an implicit discrete system can persist when the implicit function theorem no longer guarantees a smooth successor. We develop a joint Lyapunov–Schmidt reduction of the periodic equations and the determinants marking singular transitions. In unrestricted smooth one-parameter families, the generic event is a single evolution fold on a continued orbit. Simultaneous normal forms describe its interactions with saddle-node and period-doubling bifurcations, locating the singular transitions and, under backward regularity, determining the diverging multipliers. The results cover arbitrary periods and periodically non-autonomous systems. Exact realisations and applications to labour-market equilibrium, machining and implicit integration illustrate the theory.

\end{abstract}

\section{Introduction}\label{sec:intro}
Implicit discrete models describe evolution through relations between present and future states. Such relations arise, for instance, in labour-market equilibrium \cite{LazaryanLubik2019}, where future conditions enter current decisions; in machining \cite{KissBachrathy2022}, where consecutive surface errors are coupled; and in implicit integration \cite{GreenPatrickSpiteri2019}, where each new state solves a nonlinear equation. As parameters vary, a step may cease to determine a unique smooth successor even though a periodic orbit persists. A bifurcation theory must account for both the orbit and this changing evolution geometry.

We study this geometry by reducing the periodic equations together with the determinants of the next-state derivatives. Each determinant serves as a \emph{singularity marker}: its zero set locates failure of the corresponding implicit-function condition. Retaining these markers records which transitions become singular along the periodic solution set. The joint reduction remains applicable when individual steps cannot be solved smoothly for their predecessors, as the example in Section~\ref{subsec:singular-backward-example} illustrates.

For unrestricted smooth one-parameter families, we prove that a generic loss of forward evolution occurs at a nondegenerate fold of one transition (Theorem~\ref{thm:exhaustion}). Two local successor branches merge, while the periodic orbit continues uniquely. With two parameters, the joint reduction determines the geometry of collisions with a saddle-node or period doubling (Theorems~\ref{thm:ms-normal} and~\ref{thm:mixed-flip}). In each case, a curve of singular periodic orbits has quadratic contact with the periodic bifurcation curve. Cyclic equations retain the period and, in non-autonomous systems, the phase of every transition.

Our approach builds on the cyclic formulation developed by the second author \cite{Oliveira2026LS}, following the periodic equations and phase-invariance results of \cite{Oliveira18,Oliveira16}. That formulation treats bifurcations when each step determines a smooth successor. The classical reduction and singularity framework is developed in \cite{GS85}; Jepson and Spence \cite{JS89} compare equivalent reduced equations, while Vanderbauwhede \cite{Vanderbauwhede2010} adapts reduction to periodic orbits and their cyclic symmetry. The exchange symmetry used for period doubling is analysed in \cite{PeckhamKevrekidis1991}. The transition folds and genericity arguments use the singularity and transversality theory of \cite{GG73}.

Beardmore and Webster \cite{BW08} develop normal forms for singular difference-algebraic equations and obtain period doubling through doubled equations. Here the reduction transports the transition markers as well as the periodic equation, with an equivalence theory preserving their joint zero sets and recognition conditions. This gives simultaneous normal forms for the periodic and singularity loci at arbitrary periods. In the doubled problem, two separate markers identify the singular transition on either copy of the orbit; under backward regularity, their product determines the diverging multiplier.

The backward-return description connects the analysis to noninvertible maps \cite{FrouzakisKevrekidisPeckham2003} and to the forward and backward dynamics of implicit economic models \cite{GardiniEtAl2009}. Sander \cite{Sander99} develops hyperbolic sets for maps and relations, de la Llave and Lomel\'i \cite{dLL12} construct invariant manifolds for analytic difference equations, and Iglesias-Ponte et al.\ \cite{IMMP13} give a geometric formulation of implicit dynamics.

In continuous time, singularity-induced bifurcations have been studied through stability loss \cite{Venkatasubramanian1995}, matrix pencils \cite{Beardmore01}, and low-dimensional classification \cite{OvsyannikovRuan2022}; He and Barnett \cite{HeBarnett2006} examine economic models. In the generic discrete event studied here, the diverging multiplier is unstable on both sides of the singularity.

Higher-order reduced coefficients and several singularity markers provide a framework for higher codimension. Interactions with a cusp or a degenerate period doubling, and simultaneous singularity of distinct steps, are natural next problems. The codimension-one and codimension-two results proved here establish the first cases of this programme.

Section~\ref{sec:setting} sets up the periodic relations and transition geometry. Section~\ref{sec:continued} establishes continuation, spectrum and genericity. Sections~\ref{sec:mixed-saddle}--\ref{sec:mixed-spectrum} develop the mixed interactions and their spectra. Section~\ref{sec:examples} gives exact realisations at prescribed periods and applications; Section~\ref{sec:scope} identifies further directions.

\section{Periodic relations and transition folds}\label{sec:setting}
\subsection{States, phases and regularity}
Let $d,m\ge1$ be integers, and let $\Omega\subset\R^d$ and
$\Lambda\subset\R^m$ be open, with $0\in\Lambda$. The state at time
$n\in\mathbb Z$ is $x_n\in\Omega$, and $\lambda\in\Lambda$ is a
parameter held fixed along each solution. Throughout, smooth means
$C^\infty$.

An autonomous implicit system is specified by a smooth function
\begin{equation}\label{eq:auto}
 F:\Omega\times\Omega\times\Lambda\longrightarrow\R^d,
 \qquad F(x_n,x_{n+1};\lambda)=0.
\end{equation}
A pair $(x,y)$ is an \emph{admissible transition} at parameter
$\lambda$ if $F(x,y;\lambda)=0$.

For the non-autonomous case, fix a cycle of $p\ge1$ phases and
smooth functions $F_0,\ldots,F_{p-1}$ of the same type, extended
by $F_{j+p}=F_j$. The equations are
\begin{equation}\label{eq:forced}
 F_{n\bmod p}(x_n,x_{n+1};\lambda)=0.
\end{equation}
The index $n\bmod p$ is the \emph{phase}. At phase $j$,
a transition $(x,y)$ is admissible if $F_j(x,y;\lambda)=0$.
Here $p$ is the length of the prescribed phase cycle; the defining
functions may themselves repeat within that cycle. The phase labels
are retained under changes of defining equations.

A solution is a state sequence whose successive pairs are
admissible. The \emph{period} of a periodic sequence is its smallest
positive repetition length. In the autonomous case, the solution's
period is that of $(x_n)$.
In the non-autonomous case, periodicity and period refer to
$(n\bmod p,x_n)$, so the orbit period is $Np$ for an integer
$N\ge1$. A \emph{periodic orbit} is the corresponding cyclic
sequence of transitions, retaining phase in the non-autonomous case.

Set $G_j=F$ in the autonomous case and $G_j=F_{j\bmod p}$ otherwise.
Fix a reference orbit of period $\ell$, represented at $\lambda=0$
by $X_*=(x_0^*,\ldots,x_{\ell-1}^*)$, with indices read modulo
$\ell$. Thus $p\mid\ell$ in the non-autonomous case.
All local statements concern a fixed \emph{itinerary}:
neighbourhoods of these transitions, visited in their stated order
and with their phases retained. The neighbourhood of $X_*$ is chosen
disjoint from its distinct phase-compatible cyclic shifts.

Write $D_1G_j$ and $D_2G_j$ for differentiation with respect to
the present and next state, respectively. These arguments are
differentiated independently, even when the two states coincide.
At the reference orbit, put
\begin{equation}\label{eq:AB}
 A_j=D_1G_j(x_j^*,x_{j+1}^*;0),\qquad
 B_j=D_2G_j(x_j^*,x_{j+1}^*;0),\qquad 0\le j<\ell.
\end{equation}
The notation $D_y^rG_j$ denotes the $r$th derivative with respect
to the next state.

The $j$th transition is \emph{relation regular} if
$\rank[A_j\ B_j]=d$, \emph{forward regular} if $B_j$ is invertible,
and \emph{backward regular} if $A_j$ is invertible. At a fixed
parameter, relation regularity makes the nearby admissible pairs
a smooth $d$-dimensional manifold. Forward and backward regularity
allow the implicit function theorem to determine a smooth local
successor and predecessor, respectively. A transition is
\emph{forward singular} when $\det B_j=0$, so the forward
implicit-function condition is not satisfied. Norms and
orthogonality are Euclidean unless stated otherwise. For a square
matrix $T$, its corank is $\corank T=\dim\ker T$.

\subsection{Periodic equations and multipliers}\label{subsec:periodic-equations}
For $X=(x_0,\ldots,x_{\ell-1})\in\Omega^\ell$, define the map
$\Phi:\Omega^\ell\times\Lambda\to\R^{d\ell}$ by
\begin{equation}\label{eq:Phi}
 \Phi(X,\lambda)=
 \begin{pmatrix}
 G_0(x_0,x_1;\lambda)\\
 \vdots\\
 G_{\ell-1}(x_{\ell-1},x_0;\lambda)
 \end{pmatrix},\qquad x_\ell=x_0.
\end{equation}
Its zeros are solutions closing after $\ell$ steps; their periods
divide $\ell$. Thus $\ell$ is the block length in these equations,
and need not be the period of every zero.

Fix an integer $r\ge1$ and form a block of length $r\ell$ by
repeating $X_*$. Every sufficiently nearby solution of the
corresponding cyclic equations has period divisible by $\ell$.
Indeed, each cyclic shift by a number of positions not divisible
by $\ell$ changes the reference block, with phase retained when
present; these finitely many inequalities persist nearby.
In particular, the zeros near $X_*$ of the length-$\ell$ equations
have period $\ell$.

The state Jacobian $L=D_X\Phi(X_*,0)\in\R^{d\ell\times d\ell}$
acts by
\begin{equation}\label{eq:L}
 (Lv)_j=A_jv_j+B_jv_{j+1},\qquad v_\ell=v_0.
\end{equation}
In particular, $\Phi(x_0,\lambda)=G_0(x_0,x_0;\lambda)$ and $L=A_0+B_0$ when $\ell=1$. The orbit is \emph{cyclically regular} if $L$ is invertible.

At a forward-regular orbit the local successors $f_{j,\lambda}$ give
\[
 J_j=Df_{j,0}(x_j^*)=-B_j^{-1}A_j,\qquad
 P_\lambda=f_{\ell-1,\lambda}\circ\cdots\circ f_{0,\lambda},\qquad
 M=DP_0(x_0^*)=J_{\ell-1}\cdots J_0.
\]
The matrix $M$ has size $d\times d$. Its eigenvalues, counted algebraically, are the \emph{forward multipliers}. An orbit is \emph{hyperbolic} if none has modulus one; its \emph{unstable dimension} counts those with modulus greater than one \cite{Oliveira2026LS}.

\subsection{Perturbations, codimension and equivalence}\label{subsec:perturbations}

\paragraph{Allowed perturbations.}
We allow arbitrary smooth perturbations of the defining functions,
preserving autonomy or the prescribed periodic dependence on time.
An orbit is \emph{transition separated} if the ordered pairs
$(x_j^*,x_{j+1}^*)$ are all distinct in the autonomous case, or
the triples $(j\bmod p,x_j^*,x_{j+1}^*)$ are all distinct in the
non-autonomous case. An orbit of period $\ell$ can still
repeat an individual transition.

The value of a defining function and its derivatives up to a
specified finite order constitute its \emph{jet} of that order.
At separated transitions, these data can be varied independently
using perturbations with disjoint supports for each defining
function. Repeated occurrences of the same transition share the
same data and cannot be perturbed independently. Two defining
functions represent the same \emph{germ} at a transition if they
agree on some neighbourhood of that transition.

\paragraph{Codimension and genericity.}
For a fixed family, the \emph{periodic incidence set} is
$\{(X,\lambda):\Phi(X,\lambda)=0\}$. If $L$ is invertible,
it is locally the graph of $X(\lambda)$ and has dimension $m$.

The codimension counts also allow the defining functions to vary.
We consider the transition locations together with the finite jets
needed to specify the event, subject to the periodic value equations.
This larger space is the \emph{universal periodic incidence}.
Repeated transitions share a single jet, and phase is retained.
An event has codimension $r$ if it imposes $r$ additional scalar
conditions whose differentials are independent on the relevant
smooth piece, or \emph{stratum}, of this universal incidence.

Along a smooth one-parameter branch $X(\lambda)$ of periodic
solutions, a smooth scalar condition $c(X,\lambda)=0$ is crossed
\emph{transversely} at $\lambda=0$ if
\[
 c(X(0),0)=0,
 \qquad
 \left.\frac{d}{d\lambda}c(X(\lambda),\lambda)\right|_{\lambda=0}
 \ne0.
\]
Thus the condition changes to first order as the parameter varies.

A property is \emph{generic} if it holds for a \emph{residual}
set of smooth families: a set containing a countable intersection
of open dense sets in the Whitney $C^\infty$ topology.
Independent variation of the jets allows us to apply the multijet
transversality theorem \cite{GG73}. At transition-separated
configurations, a codimension-$r$ stratum has a locus of dimension
$m-r$ in a generic $m$-parameter family whenever that locus is nonempty;
it is empty if $r>m$. Section~\ref{sec:classification} treats repeated
transitions separately.

\paragraph{Allowed changes of equations and coordinates.}
A \emph{regular change of defining equations} has the form
$\widetilde G_j=\Gamma_jG_j$, where $\Gamma_j$ is a smooth
invertible $d\times d$ matrix-valued function. This leaves the
admissible pairs unchanged. These changes must respect autonomy
or the prescribed periodic dependence on time.

All coordinate changes are smooth local diffeomorphisms.
For changes of state coordinates, the same chart must be used
whenever a state is revisited; in the non-autonomous case,
this requirement applies to visits at the same phase.
When classifying the projections $(x,y)\mapsto x$ and
$(x,y)\mapsto y$ from the relation $G_j(x,y;\lambda)=0$,
we also allow independent coordinate changes on the relation
and on the target state space. These changes need not define
a conjugacy of the dynamical system.

In the reduced normal forms, new state coordinates may depend
on both states and parameters, whereas new parameters depend
only on the original parameters. A \emph{unit} is a smooth
scalar function that does not vanish near the reference point.

\subsection{Local transition geometry}\label{sec:LS}
Write $F$ for one chosen $G_j$, with $F(x_*,y_*;0)=0$, and set $A=D_xF(x_*,y_*;0)$, $B=D_yF(x_*,y_*;0)$. On the relation $\mathcal R_\lambda=\{F(x,y;\lambda)=0\}$, let $\pi_{x,\lambda}(x,y)=x$ and $\pi_{y,\lambda}(x,y)=y$. Under relation regularity, tangent vectors satisfy $A\xi+B\eta=0$, so the two projection kernels are $\{0\}\times\ker B$ and $\ker A\times\{0\}$.

\paragraph{Transition reduction.}
Suppose $\corank B=r_B\ge1$. Let $V_B,W_B\in\R^{d\times r_B}$ have orthonormal columns spanning $\ker B$ and $\ker B^T$, respectively, and put
\[
 E_B=(\ker B)^\perp,\qquad Q_B=I-W_BW_B^T,\qquad
 R_B=(B|_{E_B})^{-1}Q_B.
\]
The range equation $Q_BF(x,y_*+V_Bt+z;\lambda)=0$, with $z\in E_B$, has invertible derivative $B|_{E_B}:E_B\to\Ran B$. It determines $z=Z_B(x,t,\lambda)$, vanishing at the reference point. Define
\begin{equation}\label{eq:LS}
 Y_B=y_*+V_Bt+Z_B(x,t,\lambda),\qquad
 g_B(x,t,\lambda)=W_B^TF(x,Y_B;\lambda).
\end{equation}
Then $F=0$ is equivalent to $y=Y_B$, $g_B=0$. At the reference point,
\[
 D_tZ_B=0,\qquad D_xZ_B=-R_BA,\qquad D_tg_B=0,\qquad D_xg_B=W_B^TA.
\]
Relation regularity makes $W_B^TA$ surjective: a left null vector would annihilate both $A$ and $B$. Schur elimination in $D_yF$ gives a unit $c_B$ with
\begin{equation}\label{eq:schur}
 \det D_yF(x,Y_B;\lambda)=c_B(x,t,\lambda)\det D_tg_B(x,t,\lambda).
\end{equation}
No inverse of $A$ is used. Exchanging the state arguments gives the predecessor reduction: with orthonormal bases $V_A,W_A$ for $\ker A,\ker A^T$, solve the range equation for $Z_A\in(\ker A)^\perp$ and set
\begin{equation}\label{eq:LS-present}
 X_A=x_*+V_At+Z_A(y,t,\lambda),\qquad
 g_A(y,t,\lambda)=W_A^TF(X_A,y;\lambda).
\end{equation}
Then $F=0$ is equivalent to $x=X_A$, $g_A=0$, and $\det D_xF=c_A\det D_tg_A$ for a unit $c_A$. At the reference point, $D_tg_A=0$ and $D_yg_A=W_A^TB$ is surjective under relation regularity.

\paragraph{Fold recognition.}
For $r_B=1$, write the kernel bases as unit vectors $v_B,w_B$ and define
\begin{equation}\label{eq:a2}
 a_2=w_B^TD_y^2F(x_*,y_*;0)[v_B,v_B].
\end{equation}
An \emph{evolution fold} is a critical point of $\pi_{x,\lambda}$ with local form
\begin{equation}\label{eq:foldprojection}
 (u,s)\longmapsto(u,s^2),\qquad u\in\R^{d-1},\quad s\in\R.
\end{equation}
Here $(u,s)$ are coordinates on the relation and $(u,\eta)$ on present-state space; the \emph{critical-value boundary} is $\eta=0$. The coordinate $u$ is absent for $d=1$, and all charts may depend smoothly on parameters.

\begin{theorem}\label{thm:fold}
A relation-regular transition with $\corank B=1$ and $a_2\ne0$ is an evolution fold.
\end{theorem}
\begin{proof}
At the reference point, $(g_B)_t=0$, $(g_B)_x=w_B^TA\ne0$ and $(g_B)_{tt}=a_2$. The parameter-dependent Morse lemma gives $g_B=\varepsilon(s^2-\eta(x,\lambda))$, with $s_t\ne0$ and $\varepsilon\in\{-1,1\}$. Since $D_x\eta=-\varepsilon w_B^TA\ne0$, completing $\eta$ to present-state coordinates gives \eqref{eq:foldprojection}.
\end{proof}

There are two, one or no nearby successors according as $\eta>0$, $\eta=0$ or $\eta<0$. On the relation, \eqref{eq:schur} gives
\begin{equation}\label{eq:markerfold}
 \det D_yF=c_1(u,s,\lambda)s
\end{equation}
for a unit $c_1$, so this determinant regularly defines the critical set. A \emph{backward fold} is the corresponding fold of $\pi_{y,\lambda}$; its test is $w_A^TD_x^2F[v_A,v_A]\ne0$ for unit kernel vectors of $A,A^T$.

The fold test is representation invariant. For $\widetilde F=\Gamma F$, the extra terms in $D_y^2(\Gamma F)[v_B,v_B]$ vanish at the transition because $F=0$ and $Bv_B=0$; the corresponding left kernel vector therefore changes $a_2$ only by a nonzero factor. Coordinate invariance follows from its definition as a projection singularity.

\section{The generic singular event}\label{sec:continued}
\subsection{Continuation and crossing of a transition fold}
Let $\lambda\in\R$. If $L$ is invertible, the periodic equations determine a unique local branch $X(\lambda)$ with $X(0)=X_*$. Evaluate $A_j(\lambda),B_j(\lambda)$ along this branch and mark one index $k$. Its \emph{singularity marker} is
\begin{equation}\label{eq:sigma}
 \sigma(\lambda)=\det B_k(\lambda).
\end{equation}
Write $\adj T$ for the classical adjugate of a square matrix $T$.
Total differentiation gives the computable crossing test
\begin{align}
 X'(0)&=-L^{-1}\Phi_\lambda(X_*,0),\label{eq:branchderivative}\\
 \sigma'(0)&=\operatorname{tr}\bigl(\adj(B_k)B_k'(0)\bigr),\label{eq:sigmaprime}\\
 B_k'(0)&=D_1(D_2G_k)[x_k'(0)]+D_2(D_2G_k)[x_{k+1}'(0)]+\partial_\lambda D_2G_k.
 \label{eq:Bprime}
\end{align}
All derivatives of $G_k$ are evaluated at the reference transition. For $\corank B_k=1$, $\adj B_k$ is a nonzero multiple of $v_Bw_B^T$, so crossing is transverse exactly when $w_B^TB_k'(0)v_B\ne0$.

Collect the continued-orbit hypotheses as \textup{(E)}:
\begin{enumerate}[label=\textup{(E\arabic*)},leftmargin=4em]
 \item Every transition is relation regular; $B_k$ has corank one and all other $B_j$ are invertible.
 \item The marked fold coefficient $a_2$ in \eqref{eq:a2} is nonzero.
 \item The cyclic derivative $L$ is invertible.
 \item Along the branch supplied by (E3), $\sigma'(0)\ne0$.
\end{enumerate}
In the fold coordinates of Theorem~\ref{thm:fold}, let $s_*(\lambda)$ be the continued transition coordinate. By \eqref{eq:markerfold}, (E4) gives $s_*'(0)\ne0$, so $\nu=s_*(\lambda)$ is a local parameter. Write $\eta_*(\nu)$ for its present-state coordinate. A curve has \emph{quadratic contact} with the critical-value hypersurface $\mathcal C=\{(x,\lambda):\eta(x,\lambda)=0\}$ when a regular defining function restricts to a zero of order two on the curve.

\begin{theorem}\label{thm:contact}
Under \textup{(E)}, the periodic orbit continues uniquely and satisfies $\eta_*(\nu)=\nu^2$. At its marked present state, the nearby successors have coordinates $s=\nu$ and $s=-\nu$.
\end{theorem}
\begin{proof}
The continued transition lies on $\eta=s^2$, hence $\eta_*=\nu^2$. At the same present state, every nearby successor satisfies $s^2=\nu^2$, and the fold reconstruction gives one successor for each root. Local uniqueness of the periodic tuple follows from (E3).
\end{proof}

The transition crosses the critical set transversely, while its present state has quadratic contact with $\mathcal C$. The two successors separate to first order in $|\nu|$; only $s=\nu$ closes into the nearby periodic orbit. For nonzero small $\lambda$, the selected forward step satisfies
\begin{equation}\label{eq:step-pole}
 Df_{k,\lambda}(x_k(\lambda))
 =-\frac{\adj(B_k(0))A_k(0)}{\sigma'(0)\lambda}+O(1).
\end{equation}
Its numerator has rank one and is nonzero, since relation regularity gives $w_B^TA_k\ne0$. Thus one step derivative diverges while the periodic tuple has the finite derivative \eqref{eq:branchderivative}.

\subsection{Backward return and spectrum}\label{sec:return-spectrum}
Assume backward regularity. The local predecessor maps satisfy
\begin{equation}\label{eq:backward}
 G_j(h_{j,\lambda}(y),y;\lambda)=0.
\end{equation}
Their backward return and its derivative along a periodic branch are
\begin{equation}\label{eq:HK}
 \begin{aligned}
 H_\lambda&=h_{0,\lambda}\circ\cdots\circ h_{\ell-1,\lambda},\\
 C_j(\lambda)&=Dh_{j,\lambda}(x_{j+1}(\lambda))=-A_j(\lambda)^{-1}B_j(\lambda),\\
 K(\lambda)&=DH_\lambda(x_0(\lambda))=C_0(\lambda)\cdots C_{\ell-1}(\lambda).
 \end{aligned}
\end{equation}
After restricting the itinerary neighbourhoods, a forward block with initial state $x$ and terminal state $y$ has the complete return relation
\begin{equation}\label{eq:returnrelation}
 x=H_\lambda(y).
\end{equation}
The intermediate states are uniquely reconstructed from $y$, so periodic closure is $x=H_\lambda(x)$ even at a singular forward transition.

For $z\in\C$ and $v\in(\C^d)^\ell$, define the twisted cyclic operator by
\[
 (L_\ell(z;\lambda)v)_j=A_j(\lambda)v_j+B_j(\lambda)v_{j+1},
 \qquad v_\ell=zv_0,
\]
and put $\Delta_\lambda(z)=\det L_\ell(z;\lambda)$ and $\mathcal A(\lambda)=\prod_{j=0}^{\ell-1}\det A_j(\lambda)$. Thus $L_\ell(1;0)=L$ and $L_1(z;\lambda)=A_0(\lambda)+zB_0(\lambda)$. A \emph{fold map} has the local form \eqref{eq:foldprojection}.

\begin{proposition}\label{prop:backward-determinant}\label{prop:returnfold}
Under backward regularity,
\begin{equation}\label{eq:cyclic-backward-determinant}
 \Delta_\lambda(z)=\mathcal A(\lambda)\det(I-zK(\lambda)).
\end{equation}
If \textup{(E1)--(E2)} hold, $H_0$ is a fold map.
\end{proposition}
\begin{proof}
Factoring $A_j$ from block row $j$ leaves the operator with rows
$v_j-C_jv_{j+1}$ and closing condition $v_\ell=zv_0$.
Eliminating $v_{\ell-1},\ldots,v_1$ leaves the block
$I-zC_0\cdots C_{\ell-1}=I-zK$ on $v_0$ and identity diagonal
blocks on the other variables. These determinant-preserving row
operations give \eqref{eq:cyclic-backward-determinant} for every $z$,
including zero. This is the backward form of the twisted cyclic
identity in \cite[Theorem~2.3]{Oliveira2026LS}. It uses only
invertibility of the $A_j$ and remains valid when the blocks are
evaluated away from $\Phi=0$.

The parametrisation $y\mapsto(h_{k,0}(y),y)$ identifies the marked projection with $h_{k,0}$, a fold by Theorem~\ref{thm:fold}. Every other backward step is a local diffeomorphism, so their composition with this fold is a fold.
\end{proof}

To relate singularity markers to eigenvalues, consider a smooth
product $K(t)=\prod_{j=0}^{n-1}[-A_j(t)^{-1}B_j(t)]$, with $n\ge1$
and parameter vector $t$ near zero. Assume all $A_j$ are invertible
and $K(0)$ has a simple zero eigenvalue. Let $\kappa(t)$ continue
that eigenvalue, set $\mathcal I=\{j:\det B_j(0)=0\}$, and denote
the product of the remaining eigenvalues by $\tau_{\det}(t)$.
Define
\[
 e_{\det}(t)=\frac{(-1)^{dn}\prod_{j\notin \mathcal I}\det B_j(t)}
 {\tau_{\det}(t)\prod_{j=0}^{n-1}\det A_j(t)}.
\]

\begin{lemma}[Factorisation of the small eigenvalue]\label{lem:marker-spectral-product}
Under these assumptions, $e_{\det}$ is a smooth unit and
\begin{equation}\label{eq:marker-spectral-product}
 \kappa(t)=e_{\det}(t)\prod_{j\in\mathcal I}\det B_j(t).
\end{equation}
\end{lemma}
\begin{proof}
Division of the characteristic polynomial by its simple factor
shows that $\tau_{\det}$ is smooth and nonzero. Taking determinants
in the product for $K$ and using $\det K=\kappa\tau_{\det}$ gives
the identity. Every factor defining $e_{\det}$ is nonzero near zero.
\end{proof}

On a regular periodic orbit, only the reciprocal of $\kappa$ can
diverge near the reference point.

The cyclic determinant $\Delta_\lambda$ has degree at most $d$, and
\begin{equation}\label{eq:leading-cyclic-coefficient}
 [z^d]\Delta_\lambda(z)=(-1)^{d(\ell-1)}\prod_{j=0}^{\ell-1}\det B_j(\lambda).
\end{equation}
At $z=1$, \eqref{eq:cyclic-backward-determinant} gives $\det L=\mathcal A(0)\det(I-K(0))$. Writing $\Delta_\lambda(z)=\sum_{i=0}^d b_i(\lambda)z^i$, its degree-$d$ homogenisation is
\[
 \widehat\Delta_\lambda(Z,W)=\sum_{i=0}^d b_i(\lambda)Z^iW^{d-i}
 =\mathcal A(\lambda)\det(WI-ZK(\lambda)).
\]
In homogeneous coordinates $[Z:W]$, the affine spectral parameter
is $z=Z/W$ when $W\ne0$, and $[1:0]$ is the point at infinity.
The \emph{spectrum at infinity} is the root there, counted with
multiplicity $d-\deg\Delta_\lambda$. At a forward-regular orbit the
roots of $\Delta_\lambda$ are the forward multipliers
\cite{Oliveira2026LS}; equivalently, the local inverse returns satisfy
\begin{equation}\label{eq:inverse-return-matrices}
 M(\lambda)=K(\lambda)^{-1}.
\end{equation}

Under (E1) and backward regularity, $K(0)$ has rank $d-1$, but zero need not be algebraically simple. Impose the spectral conditions \textup{(S)}:
\begin{equation}\label{eq:spectral-clean-condition}
 \begin{gathered}
 0\text{ is an algebraically simple eigenvalue of }K(0);\\
 \text{every other eigenvalue has modulus different from one.}
 \end{gathered}
\end{equation}
These conditions define an open subset of the single-fold stratum and imply (E3). Write $\vartheta_1,\ldots,\vartheta_{d-1}$ for the nonzero eigenvalues, counted algebraically, and $\tau=\prod_{i=1}^{d-1}\vartheta_i$, with $\tau=1$ for $d=1$. Choose real null vectors $v_K,w_K$ with $w_K^Tv_K=1$, and set
\begin{equation}\label{eq:zero-eigenvalue-speed}
 a_\infty=w_K^TK'(0)v_K.
\end{equation}
Let $n_u^{\mathrm{back}}$ be the backward unstable dimension, counting the $\vartheta_i$ outside the unit circle, and let $\kappa_0(\lambda)$ continue the simple zero eigenvalue, with $\rho_\infty=\kappa_0^{-1}$ on the forward-regular branch.

\begin{theorem}\label{thm:simple-infinite-multiplier}
Under \textup{(E)}, backward regularity and \textup{(S)}, $a_\infty\ne0$ and
\begin{equation}\label{eq:multiplier-simple-pole}
 \rho_\infty(\lambda)=\frac1{a_\infty\lambda}+O(1).
\end{equation}
The remaining forward multipliers converge, with multiplicity, to $\vartheta_1^{-1},\ldots,\vartheta_{d-1}^{-1}$. The forward unstable dimension is $d-n_u^{\mathrm{back}}$ on both sides of zero.
\end{theorem}
\begin{proof}
Simple-eigenvalue differentiation gives $\kappa_0'(0)=a_\infty$ and $(\det K)'(0)=a_\infty\tau$. Since
\[
 \det K(\lambda)=(-1)^{d\ell}\frac{\prod_j\det B_j(\lambda)}{\prod_j\det A_j(\lambda)},
\]
only one numerator factor vanishes, and
\begin{equation}\label{eq:pole-coefficient-transition-data}
 a_\infty=\frac{(-1)^{d\ell}\sigma'(0)\prod_{j\ne k}\det B_j(0)}
 {\tau\prod_{j=0}^{\ell-1}\det A_j(0)}\ne0.
\end{equation}
Thus $\kappa_0(\lambda)=a_\infty\lambda+O(\lambda^2)$ gives the pole. Continuity of the other roots and \eqref{eq:inverse-return-matrices} give their reciprocal limits. The pole is unstable on both sides; of the remaining forward modes, exactly $n_u^{\mathrm{back}}$ are stable, and (S) keeps them away from the unit circle.
\end{proof}

The adjugate formula also gives the rank-one matrix pole
\begin{equation}\label{eq:return-matrix-pole}
 M(\lambda)=\frac{\adj K(0)}{(\det K)'(0)\lambda}+O(1).
\end{equation}
In terms of the cyclic coefficients, $b_d(0)=0$, $b_d'(0)\ne0$ and $b_{d-1}(0)=(-1)^{d-1}\mathcal A(0)\tau\ne0$, so
\begin{equation}\label{eq:pole-cyclic-coefficients}
 \rho_\infty(\lambda)=-\frac{b_{d-1}(0)}{b_d'(0)\lambda}+O(1).
\end{equation}
Hence $\Delta_0$ has degree $d-1$ and one simple root at infinity. Crossing zero reverses the signs of $\rho_\infty$ and $\det M$, but not the unstable dimension.

All these quantities are phase invariant: cyclic products have the same characteristic polynomial by $\det(I-zT_1T_2)=\det(I-zT_2T_1)$, and $\mathcal A$ is unchanged. The nonlinear backward returns are also locally conjugate. Across an unmarked transition use its local diffeomorphism; across the marked one use the composition of all other backward steps. In the non-autonomous case the phases are transported with the states.

\subsection{Isolation at all multiple periods}
\label{subsec:all-period-isolation}
A \emph{full orbit} of $H_\lambda$ is a sequence $(y_n)_{n\in\mathbb Z}$ with $y_{n+1}=H_\lambda(y_n)$. Condition (S) makes the continued fixed point hyperbolic also at $\lambda=0$, allowing zero eigenvalues.

\begin{proposition}\label{prop:uniform-all-period-isolation}
Under \textup{(E)}, backward regularity, and \textup{(S)}, a fixed neighbourhood contains only the continued fixed point as a full orbit of $H_\lambda$, for every sufficiently small $\lambda$.
\end{proposition}
\begin{proof}
Translate the fixed point by setting
$T_\lambda(z)=H_\lambda(x_0(\lambda)+z)-x_0(\lambda)$.
Let $E_s\oplus E_u=\R^d$ be the stable and unstable splitting of $K(0)$, and put $K_s=K(0)|_{E_s}$ and $K_u=K(0)|_{E_u}$.
Choose adapted norms and $\alpha\in(0,1)$ with $\|K_s\|\leq\alpha$ and $\|K_u^{-1}\|\leq\alpha$, omitting absent subspaces. Use the maximum norm on the sum. The remainder
$R_\lambda(z)=T_\lambda(z)-K(0)z$ satisfies, on a sufficiently small ball of radius $r>0$ and parameter interval,
\[
 \|R_\lambda(z)\|\leq\varepsilon\|z\|,
 \qquad 0<\varepsilon<1-\alpha.
\]
For a full orbit $z_n=(s_n,t_n)$ in this ball, set
$\mathcal M=\sup_{n\in\mathbb Z}\|z_n\|$.
Writing $(R_\lambda)_s,(R_\lambda)_u$ for the two components, iteration from the past in the stable equation and from the future in the unstable equation gives
\[
 s_n=\sum_{h=0}^{\infty}K_s^h(R_\lambda)_s(z_{n-1-h}),\qquad
 t_n=-\sum_{h=0}^{\infty}K_u^{-h-1}(R_\lambda)_u(z_{n+h}).
\]
The boundary terms vanish by boundedness and the contraction estimates. The series imply
$\mathcal M\leq\varepsilon\mathcal M/(1-\alpha)$, hence $\mathcal M=0$.
Finally restrict the parameter interval so that
$\|x_0(\lambda)-x_0(0)\|<r/4$.
The fixed ball of radius $r/2$ about $x_0(0)$ lies in every translated ball and contains the continued fixed point.
\end{proof}

For an integer $r\ge1$, sampling an orbit of period $r\ell$ once every $\ell$ steps and reversing the samples gives a periodic orbit of $H_\lambda$. Backward regularity reconstructs its intermediate states uniquely. The proposition therefore excludes every additional multiple-period orbit in one fixed itinerary neighbourhood, uniformly in $r$ and including $\lambda=0$.

\subsection{Genericity in one parameter}\label{sec:classification}
We classify singular events at each fixed period in the unrestricted class of Section~\ref{sec:setting}. A \emph{type} records the transition projection singularity and its position in the universal periodic incidence of Section~\ref{subsec:perturbations}; spectral data distinguish open subsets.

\paragraph{Repeated transitions.}
Fix an equality pattern among the states. In the non-autonomous case equality labels include phase; different phases retain separate coordinates even when their numerical states coincide. Let $n_{\mathrm v}$ be the number of distinct labels and $n_{\mathrm e}$ the number of distinct transitions. Form the directed graph with these labels as vertices and transitions as edges. The periodic word is a closed walk visiting every vertex and edge. Independent jets are available at distinct edges, but not at repeated uses of one edge.

\begin{lemma}\label{lem:repeated-transitions}
At each fixed period, a generic one-parameter family has no orbit containing both a singular next-state derivative and a repeated transition.
\end{lemma}
\begin{proof}
The graph of the cyclic word is strongly connected: following the word gives a directed path between any two visited vertices. Every vertex has at least one outgoing edge, so $n_{\mathrm e}\ge n_{\mathrm v}$. If equality holds, every vertex has exactly one outgoing edge, and strong connectivity makes the graph a single directed cycle. The word represents one period of the orbit, so it traverses that cycle once. An orbit with a repeated transition must therefore satisfy $n_{\mathrm e}\ge n_{\mathrm v}+1$.

For its fixed equality pattern, the distinct state coordinates range over an open subset of a product of $n_{\mathrm v}$ copies of $\Omega$. This gives dimension $dn_{\mathrm v}$. There is one vector orbit equation for each of the $n_{\mathrm e}$ distinct transition germs, hence $dn_{\mathrm e}$ scalar value conditions. Repeated uses of the same germ do not give additional independent equations. Values and first derivatives at distinct germs can be perturbed independently by smooth functions supported in disjoint neighbourhoods, separately for each phase.

At one marked germ, suppose its next-state derivative has corank exactly $r\ge1$. This is a condition of codimension $r^2$ on that matrix: choose an invertible $(d-r)\times(d-r)$ block and require the complementary $r\times r$ Schur complement to vanish. These are $r^2$ independent scalar equations, additional to the value equations.

For $m$ parameters, apply multijet transversality to the submanifold of the multijet bundle defined by the value and corank conditions, together with the source-coordinate constraints identifying successive endpoints and the common parameter. These identifications leave $dn_{\mathrm v}+m$ free source coordinates; the jets at distinct transition germs remain independently variable. For a residual set of defining families, the resulting stratum therefore has expected dimension
\[
 dn_{\mathrm v}+m-dn_{\mathrm e}-r^2
   =m-d(n_{\mathrm e}-n_{\mathrm v})-r^2.
\]
For $m=1$ and an orbit containing a repeated transition, this dimension is at most $-d$, so the transverse stratum is empty.

There are finitely many equality patterns, marked edges, and coranks at a fixed orbit length. Apply transversality on a countable cover by relatively compact coordinate neighbourhoods and intersect the resulting residual sets. This excludes all these patterns simultaneously. The argument for non-autonomous equations uses independent perturbations within each phase and is unchanged.
\end{proof}

\paragraph{Classification.}
Call an event \emph{elementary} if it satisfies (E), every $A_j$ is invertible, and the spectral conditions (S) hold. Condition (S) already implies (E3) by Proposition~\ref{prop:returnfold}.

We use semialgebraic projection and smooth stratification \cite[Chs.~2 and~9]{BCR98} to handle the exceptional spectral sets.

\begin{theorem}\label{thm:exhaustion}
At each fixed period, every loss of forward evolution in a generic one-parameter family of the unrestricted class is elementary.
\end{theorem}
\begin{proof}
Lemma~\ref{lem:repeated-transitions} removes repeated transitions. Work at separated configurations, with the transition locations and the finite jets of the defining functions as variables. Independent variation of their values shows that $\Phi=0$ imposes exactly $d\ell$ scalar conditions. Once these value equations are imposed, the first and higher derivatives at the distinct transition points remain free jet data.

The corank-$r$ matrix stratum has codimension $r^2$, as in the lemma. Thus a singular next-state derivative imposes one condition precisely on its corank-one stratum; corank at least two imposes at least four. Two distinct singular transition germs impose two independent determinant conditions. A generic one-parameter singular event consequently has a single marked matrix $B_k$ of corank one, with every other $B_j$ invertible.

On this corank-one stratum, choose smoothly varying unit null vectors $v_B,w_B$ locally. With its first jet fixed, a perturbation of the second derivative in the direction
\[
 (\xi,\zeta)\longmapsto w_B(v_B^T\xi)(v_B^T\zeta),
 \qquad \xi,\zeta\in\R^d,
\]
changes $a_2=w_B^TD_y^2G_k[v_B,v_B]$ freely. Its vanishing therefore imposes one more scalar condition. Likewise, singularity of any $A_j$ is an additional condition on independent first-derivative data. Outside these exceptional strata, all $A_j$ are invertible; this also ensures relation regularity. The remaining state conditions are precisely the fold assumptions (E1)--(E2) and backward regularity.

On this open part of the corank-one stratum, $K=C_0\cdots C_{\ell-1}$ has exactly one rank-$(d-1)$ factor and otherwise invertible factors. Varying $B_k$ within its corank-one stratum, while holding the other first derivatives fixed, makes $K$ vary freely over the rank-$(d-1)$ matrix stratum. Indeed, the variation is multiplication of $B_k$ on the left and right by fixed invertible matrices.

Within this matrix stratum, algebraic multiplicity greater than one for the zero eigenvalue is a proper algebraic condition: with $t$ denoting the polynomial indeterminate, the coefficient of its linear term in $\det(tI-K)$ must vanish. It is not identically zero, as seen from a diagonal matrix with one zero and all other entries nonzero. Thus this failure has positive codimension within the stratum.

The presence of an eigenvalue on the unit circle is also exceptional there. The subset is semialgebraic: write a complex root as $u+iv$, impose $u^2+v^2=1$, and set the real and imaginary parts of the characteristic polynomial equal to zero; elimination of $u,v$ preserves semialgebraicity. It has empty interior because multiplying any rank-$(d-1)$ matrix by an arbitrarily close positive scalar can move all its nonzero eigenvalues away from the unit circle. There are only finitely many prohibited scaling factors. A semialgebraic subset with empty interior admits a stratification by smooth pieces of positive codimension. Hence this spectral failure imposes at least one more condition. In particular, the eigenvalue $1$, and therefore singularity of $L$, is not forced by the forward singularity. When $d=1$, $K(0)=0$ and neither spectral exception occurs.

All these state and spectral exceptions have codimension at least two in the universal periodic incidence. For the remaining fold stratum, transversality means that the augmented map $(\Phi,\det B_k)$ has surjective derivative with respect to $(X,\lambda)$. Since $D_X\Phi=L$ is invertible, elimination of the $X$ columns leaves the scalar derivative $\sigma'(0)$ along the continued branch. Transversality is therefore exactly (E4).

The required parameter derivative is free jet data. For example, choose a cut-off function equal to one near the marked transition and multiply it by
\[
 \lambda\,w_Bv_B^T(y-x_{k+1}^*).
\]
Adding a scalar multiple of this expression to $G_k$ preserves its value, first state jet, and partial parameter value at the reference point. It therefore preserves $X'(0)$ in \eqref{eq:branchderivative}, while changing $w_B^TB_k'(0)v_B$ freely. At separated germs, the support does not affect any other evaluated transition.

Apply multijet transversality \cite[Ch.~II, \S4]{GG73} to these
smooth strata, with the common-parameter and endpoint identifications
used in Lemma~\ref{lem:repeated-transitions}. The separated
configurations have $d\ell+1$ free source coordinates. The
$d\ell$ value equations and any exceptional stratum of relative
codimension at least two therefore give negative expected dimension;
the single-fold stratum has expected dimension zero. Thus generic
families avoid the exceptional strata and meet the fold stratum
transversely. A countable cover by relatively compact coordinate
neighbourhoods and intersection of the resulting residual sets give
the assertion at the prescribed period.
\end{proof}

The marked position and phase label the same codimension-one type. Its finite forward spectrum may have any number from zero to $d-1$ of stable modes; the escaping multiplier is unstable on both sides. Table~\ref{tab:codimension} collects the exceptional codimension bounds.

\begin{table}[htbp]
\centering
\small
\renewcommand{\arraystretch}{1.08}
\begin{tabular}{@{}>{\raggedright\arraybackslash}p{0.57\textwidth}@{\hspace{1em}}>{\raggedright\arraybackslash}p{\dimexpr0.43\textwidth-1em\relax}@{}}
\toprule
Condition at a periodic orbit & Codimension\\
\midrule
One forward fold transition & One; the case classified here\\
Vanishing fold quadratic coefficient & At least two, including the forward singularity\\
Two distinct singular transition germs & At least two\\
A forward singularity and a backward singularity & At least two\\
A forward singularity and singular cyclic compatibility & At least two\\
A multiple zero eigenvalue of the backward return & At least two, on the single-singularity locus\\
A forward singularity and a unit-modulus finite mode & At least two\\
Corank at least two in one $B_j$ & At least four\\
\bottomrule
\end{tabular}
\caption{Codimension bounds in the universal periodic incidence at transition-separated configurations.}
\label{tab:codimension}
\end{table}

\paragraph{Periods, phases and time reversal.}\label{sec:autonomous}\label{sec:nonautonomous}
The theorem applies to autonomous and non-autonomous orbits, using independent jets within each phase. Cyclic relabelling preserves codimension. The time-reversed relations
\[
 F^{\rm rev}(x,y;\lambda)=F(y,x;\lambda),\qquad
 F_j^{\rm rev}(x,y;\lambda)=F_{(-j-1)\bmod p}(y,x;\lambda)
\]
generate $\widetilde x_n=x_{-n}$ and exchange $A$ and $B$, giving the corresponding classification of backward singularities.

\section{An evolution fold meeting a saddle-node}
\label{sec:mixed-saddle}

Let $\lambda=(\lambda_1,\lambda_2)$ and fix a reference orbit $X_*$. We study the joint degeneracy of the map $\Phi$ and one transition. The operator reduced here is $L=D_X\Phi(X_*,0)$; the transition reduction of Section~\ref{sec:LS} used $B_k$.

\subsection{A cyclic equation with a singularity marker}

Suppose $\corank L=1$. Choose unit vectors $v_L,w_L$ spanning $\ker L,\ker L^T$, and define
\[
 E_L=(\ker L)^\perp,\qquad Q_L=I-w_Lw_L^T,
 \qquad R_L=(L|_{E_L})^{-1}Q_L:\R^{d\ell}\to E_L.
\]
The cyclic range equation $Q_L\Phi(X_*+sv_L+z,\lambda)=0$, $z\in E_L$, determines $z=Z_L(s,\lambda)$ because its $z$-derivative is $L|_{E_L}:E_L\to\Ran L$ \cite{Oliveira2026LS}. Define

\begin{equation}\label{eq:ms-reduction}
 \mathcal X(s,\lambda)=X_*+sv_L+Z_L(s,\lambda),\qquad
 h(s,\lambda)=w_L^T\Phi(\mathcal X(s,\lambda),\lambda).
\end{equation}
The nearby periodic solutions are exactly $\mathcal X(s,\lambda)$ with $h=0$. Write $\mathcal X_j$ for its $j$th state and $\mathcal X_\ell=\mathcal X_0$. At the origin, $(Z_L)_s=0$, $\mathcal X_s=v_L$, and $h=h_s=0$.

For a marked transition $k$, define on the full range reconstruction

\begin{equation}\label{eq:ms-marker}
 \sigma(s,\lambda)=
 \det D_2G_k(\mathcal X_k(s,\lambda),
             \mathcal X_{k+1}(s,\lambda);\lambda).
\end{equation}
On $h=0$, the equation $\sigma=0$ identifies the singular periodic transitions. Assume the marked transition is a relation-regular evolution fold, all other $B_j$ are invertible, and

\begin{equation}\label{eq:ms-conditions}
 h_{ss}\ne0,\qquad \sigma_s\ne0,\qquad
 \rank D_{(s,\lambda)}(h,h_s,\sigma)=3.
\end{equation}
All derivatives are evaluated at $(s,\lambda)=(0,0)$, where $h=h_s=\sigma=0$. The inequalities express quadratic turning of the periodic equation and transverse crossing of the marker along its kernel; the rank condition gives an independent unfolding. Denote these assumptions, together with $\corank L=1$, by \textup{(SN)}.

We use local changes of state and parameter coordinates of the form $(s,\lambda)\mapsto(u(s,\lambda),\psi(\lambda))$, with $u_s(0,0)\ne0$ and $\psi$ a local parameter diffeomorphism fixing the origin. Write $(\mu,\nu)=\psi(\lambda)$ and compose functions with the inverse change. The \emph{periodic-fold locus} is the parameter projection of $h=h_s=0$; the \emph{singular-orbit locus} is the projection of $h=\sigma=0$. Quadratic contact means that a regular defining function of either curve vanishes to order two on a regular parametrisation of the other. Below, $\alpha_h,\alpha_\sigma$ denote units.

\begin{theorem}\label{thm:ms-normal}
Under \textup{(SN)}, a coordinate change of this form gives
\begin{equation}\label{eq:ms-normal}
 h=\alpha_h(u,\mu,\nu)(u^2-\mu),\qquad
 \sigma=\alpha_\sigma(u,\mu,\nu)(u-\nu).
\end{equation}
The periodic-fold and singular-orbit loci are $\mu=0$ and $\mu=\nu^2$, with quadratic contact.
\end{theorem}

\begin{proof}
The parameter-dependent Morse lemma puts the scalar equation into saddle-node form. Explicitly, let $c(\lambda)$ solve $h_s(c(\lambda),\lambda)=0$. Taylor's formula gives
\[
 h(s,\lambda)=h(c(\lambda),\lambda)
                +(s-c(\lambda))^2b_h(s,\lambda),
 \qquad b_h(0,0)=\tfrac12h_{ss}(0,0)\ne0.
\]
With $\epsilon$ the sign of $b_h$, set
$u=(s-c(\lambda))\sqrt{|b_h(s,\lambda)|}$ and
$\mu(\lambda)=-\epsilon h(c(\lambda),\lambda)$.
Then $h=\epsilon(u^2-\mu(\lambda))$ and $u_s(0,0)\ne0$.

In $(u,\lambda)$ coordinates, $\sigma_u\ne0$, so $\sigma=0$ has a unique root $u=\nu_0(\lambda)$. Division gives
\[
 \sigma=(u-\nu_0(\lambda))b_\sigma(u,\lambda),\qquad
 b_\sigma=\int_0^1\sigma_u\bigl(\nu_0(\lambda)+t(u-\nu_0(\lambda)),\lambda\bigr)\,\dd t,
\]
with $b_\sigma$ a unit. At the origin, $\dd h$, $\dd h_s$ and $\dd\sigma$ are nonzero multiples of $\dd\mu$, $\dd u$ and $\dd u-\dd\nu_0$, respectively; the original derivative $h_s$ contributes the nonzero factor $u_s$. Their independence is equivalent to independence of $\dd\mu,\dd\nu_0$ on parameter space. Thus $\lambda\mapsto(\mu,\nu_0)$ is a local diffeomorphism. Taking $\nu=\nu_0$ proves \eqref{eq:ms-normal}. Finally, $h=h_s=0$ gives $u=\mu=0$, whereas $h=\sigma=0$ gives $u=\nu$, $\mu=\nu^2$; the contact is quadratic.
\end{proof}

For $\mu>0$ the periodic roots are $u=\pm\sqrt\mu$; for $\mu<0$ there are none. On the punctured parabola only $u=\nu$ is singular. These roots represent distinct orbits of period $\ell$, by the local conventions of Section~\ref{sec:setting}.

In this two-parameter family, the rank condition makes the periodic incidence $h=0$ a smooth surface. On it, $h_s=0$ and $\sigma=0$ are independent conditions, so their common zero is an isolated event of codimension two.

\subsection{Recognition coefficients}
For $i\in\{1,2\}$, let $e_i$ be the $i$th coordinate vector of $\R^2$ and put $U_i=-R_L\Phi_{\lambda_i}(X_*,0)$. The coefficient formulas of \cite[Propositions~5.1--5.2]{Oliveira2026LS}, in these coordinates, are
\begin{align}
 h_{ss}&=w_L^TD_X^2\Phi[v_L,v_L],&
 h_{\lambda_i}&=w_L^T\Phi_{\lambda_i},\label{eq:ms-derivatives}\\
 h_{s\lambda_i}&=w_L^T\bigl(D_X^2\Phi[v_L,U_i]+(D_X\Phi_{\lambda_i})v_L\bigr).\nonumber
\end{align}
Only the cyclic range inverse $R_L$ enters their derivation, so the formulas remain valid for singular $B_k$.

For the marker calculation define the matrix-valued function
$\mathcal B_k(x,y,\lambda)=D_2G_k(x,y;\lambda)$.
If $\xi=(\xi_0,\ldots,\xi_{\ell-1})\in\R^{d\ell}$ and
$\eta=(\eta_1,\eta_2)\in\R^2$, its total variation along these state
and parameter directions is the $d\times d$ matrix
\[
 T_k[\xi,\eta]
 =D_x\mathcal B_k[\xi_k]
       +D_y\mathcal B_k[\xi_{k+1}]
       +\sum_{i=1}^2\eta_i\partial_{\lambda_i}\mathcal B_k.
\]
The derivatives on the right are evaluated at
$(x_k^*,x_{k+1}^*,0)$, and $\xi_\ell=\xi_0$.
The adjugate formula for the differential of a determinant now gives
\begin{equation}\label{eq:ms-marker-derivatives}
 \sigma_s=\operatorname{tr}\bigl(\adj(B_k)T_k[v_L,0]\bigr),\qquad
 \sigma_{\lambda_i}
 =\operatorname{tr}\bigl(\adj(B_k)T_k[U_i,e_i]\bigr).
\end{equation}
The rank condition in \eqref{eq:ms-conditions} is the explicit test
\begin{equation}\label{eq:ms-recognition-matrix}
 \det\begin{pmatrix}
 0&h_{\lambda_1}&h_{\lambda_2}\\
 h_{ss}&h_{s\lambda_1}&h_{s\lambda_2}\\
 \sigma_s&\sigma_{\lambda_1}&\sigma_{\lambda_2}
 \end{pmatrix}\ne0.
\end{equation}

\subsection{Invariance of the marked reduction}

For two marked reductions $(h,\sigma)$ and $(\widetilde h,\widetilde\sigma)$, let $\chi(s,\lambda)=(\theta(s,\lambda),\psi(\lambda))$ be a local diffeomorphism of this form, $\alpha,\beta$ units and $\gamma$ smooth. The relevant equivalence is

\begin{equation}\label{eq:ms-equivalence}
 \widetilde h\circ\chi=\alpha h,\qquad
 \widetilde\sigma\circ\chi=\beta\sigma+\gamma h.
\end{equation}
The term $\gamma h$ permits different extensions of the marker off the periodic surface.

\begin{proposition}\label{prop:ms-invariance}
Conditions \eqref{eq:ms-conditions} and the quadratic contact are
invariant under \eqref{eq:ms-equivalence}. Changes of reduction
data, regular defining equations, state coordinates, and starting
phase induce this equivalence.
\end{proposition}

\begin{proof}
First keep the reduced coordinates fixed. At
$h=h_s=\sigma=0$, multiplication by the units in
\eqref{eq:ms-equivalence} preserves the nonvanishing of
$h_{ss}$ and $\sigma_s$. The relevant differentials satisfy
\[
 \begin{aligned}
 \dd(\alpha h)&=\alpha\,\dd h,\\
 \dd((\alpha h)_s)&=\alpha\,\dd h_s+\alpha_s\,\dd h,\\
 \dd(\beta\sigma+\gamma h)
       &=\beta\,\dd\sigma+\gamma\,\dd h.
 \end{aligned}
\]
This is an invertible triangular transformation of the three
covectors. A coordinate change of this form adds nonzero
chain-rule factors and multiples of these rows. Thus it also
preserves the two inequalities and the rank condition. The
equivalence identifies the periodic and singular subsets and
applies the same parameter diffeomorphism to their projections;
their order of contact is consequently unchanged.

For two choices of reduction data, evaluate the second kernel coordinate on the first reconstruction. Its nonzero kernel derivative defines the scalar coordinate change. Uniqueness in the range solve identifies the reconstructed solutions on $h=0$. Since $\dd h\ne0$, both scalar equations are submersions with the same zero set, and division in a coordinate transverse to that set gives $\widetilde h\circ\chi=\alpha h$. The markers agree there, so their difference is $\gamma h$. This is the scalar reduction comparison used in \cite{Oliveira2026LS}, now retaining the transition marker.

For $\widetilde G_j=\Gamma_jG_j$, with $\Gamma_j$ smooth and invertible, the map $\Phi$ is multiplied by an invertible block matrix. On $G_j=0$,
\[
 \det D_2\widetilde G_j=(\det\Gamma_j)\det D_2G_j.
\]
Thus the scalar reduction comparison applies, and the displayed identity supplies $\beta$ on the periodic surface. Extending it smoothly gives the remaining multiple of $h$. State-coordinate changes multiply the markers by coordinate determinant units. Cyclic relabelling permutes the equations and transports the marked transition; in the non-autonomous case it also shifts the sequence of defining functions. Each change therefore induces \eqref{eq:ms-equivalence}.
\end{proof}

\subsection{The finite central multiplier and the escaping multiplier}

Assume backward regularity and that $K_*=DH_0(x_0^*)$ has simple eigenvalues $0,1$, with all others nonzero and off the unit circle; in particular $d\ge2$. On the periodic surface define $K(s,\lambda)=DH_\lambda(\mathcal X_0(s,\lambda))$. Let $\kappa_0,\eta_c$ continue these two eigenvalues, and set $\rho_\infty=\kappa_0^{-1}$, $m_c=\eta_c^{-1}$ on the forward-regular part. The surface $\mu=u^2$ has coordinates $(u,\nu)$; $c_0,c_c$ below are units on it.

\begin{proposition}\label{prop:ms-spectrum}
Under \textup{(SN)} and the spectral assumptions above, there
are units $c_0,c_c$ such that, on the forward-regular part of
$\mu=u^2$,
\begin{equation}\label{eq:ms-spectrum}
 \rho_\infty=\frac{1}{c_0(u,\nu)(u-\nu)},\qquad
 m_c-1=c_c(u,\nu)u.
\end{equation}
Every sufficiently nearby forward-regular periodic orbit is
unstable.
\end{proposition}

\begin{proof}
Lemma~\ref{lem:marker-spectral-product} and \eqref{eq:ms-normal} give $\kappa_0=c_0(u,\nu)(u-\nu)$, proving the pole law by reciprocation.

Schur elimination in $D_X\Phi$ gives $\det D_X\Phi=h_s$ times a unit. On $h=0$, \eqref{eq:ms-normal} makes $h_s=u$ times a unit, including the nonzero factor $u_s$. By \eqref{eq:cyclic-backward-determinant}, $\det D_X\Phi=\mathcal A\det(I-K)$, and only $1-\eta_c$ can vanish. Thus $1-\eta_c$, and hence $m_c-1$, equals $u$ times a unit. Finally, $|\kappa_0|<1$ near the reference orbit, so every regular orbit has $|\rho_\infty|>1$.
\end{proof}

The central mode lies on opposite sides of $1$ on the two roots $u=\pm\sqrt\mu$, while the pole mode makes both full orbits unstable.

\section{An evolution fold and period doubling}
\label{sec:mixed-flip}

Let $\lambda\in\R^2$ and fix a reference orbit with one marked evolution fold at index $k$, all other $B_j$ being invertible. The doubled equations \cite{PeckhamKevrekidis1991,Oliveira2026LS} repeat this transition twice, and we retain both determinant markers. With the twisted operators of Section~\ref{sec:return-spectrum}, set $L_+=L_\ell(1;0)$ and $L_-=L_\ell(-1;0)$ at the reference orbit and assume

\begin{equation}\label{eq:mf-linear-conditions}
 \det L_+\ne0,\qquad \corank L_-=1.
\end{equation}
The first condition gives a unique \emph{primary branch} $X^{\mathrm{pr}}(\lambda)=(x_0^{\mathrm{pr}},\ldots,x_{\ell-1}^{\mathrm{pr}})$ with $X^{\mathrm{pr}}(0)=X_*$. The second singles out one anti-periodic direction, with boundary condition $v_\ell=-v_0$. Neither requires a forward return at the reference parameter.

\subsection{Reduction of the doubled cyclic equations}

Let $R=(R_0,\ldots,R_{\ell-1})$ and
$S=(S_0,\ldots,S_{\ell-1})$ be two state blocks near $X_*$; each
component $R_j$ and $S_j$ belongs to $\R^d$. Define
$\mathcal D(R,S;\lambda)\in\R^{d\ell}$ by
\begin{equation}\label{eq:mf-block-residual}
 \mathcal D_j(R,S;\lambda)=
 \begin{cases}
 G_j(R_j,R_{j+1};\lambda),&0\leq j<\ell-1,\\
 G_{\ell-1}(R_{\ell-1},S_0;\lambda),&j=\ell-1.
 \end{cases}
\end{equation}
For $\ell=1$ only the closing transition is present. Define the doubled map by
\[
 \Psi((R,S),\lambda)
   =\bigl(\mathcal D(R,S;\lambda),\mathcal D(S,R;\lambda)\bigr)
   \in\R^{2d\ell}.
\]
Its zeros are labelled solutions closing after $2\ell$ steps, including repetitions of the primary orbit. The block exchange $\mathcal S(R,S)=(S,R)$ satisfies $\Psi\circ\mathcal S=\mathcal S\circ\Psi$. In the non-autonomous case this follows from $p\mid\ell$. At the repeated reference orbit,
\[
 D_{(R,S)}\Psi[\zeta,\zeta]=(L_+\zeta,L_+\zeta),\qquad
 D_{(R,S)}\Psi[\zeta,-\zeta]=(L_-\zeta,-L_-\zeta),\quad
 \zeta\in\R^{d\ell}.
\]
The splitting uses only the cyclic equations and remains valid with singular $B_k$.

Choose Euclidean unit vectors $v_-\in\ker L_-$ and
$w_-\in\ker L_-^T$. Write $v_-=(v_{-,0},\ldots,v_{-,\ell-1})$, with
$v_{-,j}\in\R^d$, and define
\[
 E_- = v_-^\perp\subset\R^{d\ell},\qquad Q_-=I-w_-w_-^T.
\]
The restriction $L_-|_{E_-}:E_-\to\Ran L_-$ is invertible.

Write the mean displacement as $v\in\R^{d\ell}$ and the half-difference as $w$:
\[
 R=X^{\mathrm{pr}}(\lambda)+v+w,
 \qquad S=X^{\mathrm{pr}}(\lambda)+v-w,
 \qquad w=sv_-+z,\quad s\in\R,\quad z\in E_-.
\]
The scalar $s$ is the anti-periodic amplitude. Define the symmetric and antisymmetric maps by
\begin{align*}
 \mathcal E(v,w,\lambda)&=\tfrac12\{
       \mathcal D(R,S;\lambda)+\mathcal D(S,R;\lambda)\},\\
 \mathcal O(v,w,\lambda)&=\tfrac12\{
       \mathcal D(R,S;\lambda)-\mathcal D(S,R;\lambda)\}.
\end{align*}
The first map is even in $w$ and the second is odd. The range equations are
\begin{equation}\label{eq:mf-range}
 \mathcal E(v,sv_-+z,\lambda)=0,\qquad
 Q_-\mathcal O(v,sv_-+z,\lambda)=0.
\end{equation}
Their derivative in $(v,z)$ at the origin is the isomorphism
$\diag(L_+,L_-|_{E_-}):\R^{d\ell}\times E_-\to\R^{d\ell}\times\Ran L_-$.
Thus they determine unique smooth $v(s,\lambda)$ and $z(s,\lambda)$, and uniqueness gives
\begin{equation}\label{eq:mf-parity}
 v(-s,\lambda)=v(s,\lambda),\qquad z(-s,\lambda)=-z(s,\lambda),
 \qquad v(0,\lambda)=z(0,\lambda)=0.
\end{equation}

Substituting these functions in the expressions for $R$ and $S$
defines the reconstruction
$Y(s,\lambda)=(R(s,\lambda),S(s,\lambda))\in\R^{2d\ell}$.
Its components are ordered as
$Y_j=R_j$ and $Y_{j+\ell}=S_j$ for $0\leq j<\ell$; their indices
are read modulo $2\ell$. The scalar function left after solving
\eqref{eq:mf-range} is
\begin{equation}\label{eq:mf-reduced}
 g(s,\lambda)
   =w_-^T\mathcal O(v(s,\lambda),sv_-+z(s,\lambda),\lambda).
\end{equation}
In particular,
\[
 \Psi(Y(s,\lambda),\lambda)=(w_-\,g(s,\lambda),-w_-\,g(s,\lambda)).
\]
Thus $\Psi=0$ is equivalent to $g=0$, with amplitude recovered as $s=v_-^T(R-S)/2$. Block exchange gives $Y(-s,\lambda)=\mathcal S Y(s,\lambda)$ and makes $g$ odd.

Whitney's lemma for smooth even functions, with parameters \cite{Whitney1943}, gives a smooth function $\omega(r,\lambda)$ near $(0,0)$ such that
\begin{equation}\label{eq:mf-odd-factorisation}
 g(s,\lambda)=s\,\omega(s^2,\lambda),\qquad \omega(0,0)=0.
\end{equation}
Here $r=s^2$, with $\omega$ smoothly extended to negative $r$. Indeed, division of the odd function $g$ by $s$ gives a smooth even function, and $\omega(0,0)=w_-^TL_-v_-=0$.

\subsection{Two markers and the normal form}
Put $a_3=\omega_r(0,0)=g_{sss}(0,0)/6$.

Define two scalar determinant functions on the reconstruction by
\begin{align}
 \sigma_+(s,\lambda)
 &=\det D_2G_k(Y_k(s,\lambda),Y_{k+1}(s,\lambda);\lambda),
 \label{eq:mf-marker-plus}\\
 \sigma_-(s,\lambda)
 &=\det D_2G_k(Y_{k+\ell}(s,\lambda),Y_{k+\ell+1}(s,\lambda);\lambda)
   =\sigma_+(-s,\lambda).
 \label{eq:mf-marker-minus}
\end{align}
Each function detects loss of the next-state invertibility condition
at one copy of the distinguished transition. They are defined even
when the reconstruction does not satisfy the final scalar equation
$g=0$. On the primary branch they agree; write their common value as
$\sigma_{\mathrm{pr}}(\lambda)=\sigma_+(0,\lambda)$.

Set $c_{\mathrm{PD}}=\partial_s\sigma_+(0,0)$. The mixed nondegeneracy conditions are
\begin{equation}\label{eq:mf-nondegeneracy}
 a_3\ne0,\qquad c_{\mathrm{PD}}\ne0,\qquad
 \rank D_\lambda\bigl(\omega(0,\lambda),\sigma_{\mathrm{pr}}(\lambda)\bigr)_{\lambda=0}=2.
\end{equation}
They require a nonzero cubic term, first-order separation of the markers and independent parameter unfolding. In particular, $c_{\mathrm{PD}}\ne0$ is additional to the evolution-fold condition.

Denote the reference-orbit and fold assumptions, \eqref{eq:mf-linear-conditions} and \eqref{eq:mf-nondegeneracy} by \textup{(PD)}.

Define $\mu(\lambda)$ and
$\nu(\lambda)$ to be the unique small roots of
\[
 \omega(\mu(\lambda),\lambda)=0,
 \qquad \sigma_+(\nu(\lambda),\lambda)=0.
\]
They exist and are smooth by $a_3\ne0$ and $c_{\mathrm{PD}}\ne0$, respectively.
The first is a root in the squared-amplitude variable $r$; the second
is a root in the signed amplitude $s$. Both vanish at $\lambda=0$.

The functions $u_g,u_\sigma$ below are smooth units.
\begin{theorem}\label{thm:mixed-flip}
Under \textup{(PD)}, $(\mu,\nu)$ are local parameter coordinates.
The nonprimary orbits have period $2\ell$ and satisfy
$s^2=\mu>0$. They contain a singular transition exactly on
$\mu=\nu^2>0$.
The simultaneous normal form is
\begin{align}
 g(s,\mu,\nu)&=su_g(s^2,\mu,\nu)(s^2-\mu),
 \label{eq:mf-normal-periodic}\\
 \sigma_+(s,\mu,\nu)&=u_\sigma(s,\mu,\nu)(s-\nu),
 \label{eq:mf-normal-plus}\\
 \sigma_-(s,\mu,\nu)&=u_\sigma(-s,\mu,\nu)(-s-\nu).
 \label{eq:mf-normal-minus}
\end{align}
\end{theorem}

\begin{proof}
Implicit differentiation gives
\[
 D_\lambda\mu(0)=-a_3^{-1}D_\lambda\omega(0,0),\qquad
 D_\lambda\nu(0)=-c_{\mathrm{PD}}^{-1}D_\lambda\sigma_{\mathrm{pr}}(0).
\]
The rank condition therefore makes $\lambda\mapsto(\mu,\nu)$
a local diffeomorphism. From this point on, a function originally
written with parameter $\lambda$ is written with $(\mu,\nu)$
by composition with the inverse parameter change. In particular,
$Y(s,\mu,\nu)$ denotes the same reconstruction in the new parameters.

Division by the simple roots gives the displayed factorisations. For example,

\[
 \omega(r,\lambda)=(r-\mu(\lambda))
 \int_0^1\omega_r(\mu(\lambda)+t(r-\mu(\lambda)),\lambda)\,\dd t.
\]
The integral defines $u_g$ after the parameter change, and its value
at the origin is $a_3$. The same integral division applied to
$\sigma_+$ defines $u_\sigma$, with $u_\sigma(0,0,0)=c_{\mathrm{PD}}$. Both are therefore
nonzero nearby.

Besides $s=0$, which reconstructs the repeated primary orbit, the
only roots are $s=\pm\sqrt\mu$ for $\mu>0$. These two labelled
solutions represent one orbit with its two blocks exchanged.
Neither has period $\ell$, since $R=S$ would force
$s=v_-^T(R-S)/2=0$.

By the local period argument of Section~\ref{subsec:periodic-equations},
a nearby solution closing after $2\ell$ steps has period $\ell$ or
$2\ell$. The nonzero roots therefore have period $2\ell$.

On these roots, the marker product is
\begin{equation}\label{eq:mf-marker-product}
 \sigma_+\sigma_-
   =u_\sigma(s,\mu,\nu)u_\sigma(-s,\mu,\nu)(\nu^2-\mu).
\end{equation}
Every other transition remains forward regular. Thus the doubled
orbit is singular precisely on $\mu=\nu^2>0$. Exactly one of its
two marked transitions is singular there: simultaneous vanishing
of the two markers would give $s=\nu=0$. A fixed nonzero minor
preserves corank one near the reference singular transition, and
the nonzero fold coefficient also persists. Its singularity is
therefore still an evolution fold.
\end{proof}

The primary periodic-degeneracy and singularity curves are $\mu=0$
and $\nu=0$, respectively, since $g_s(0,\mu,\nu)=-\mu u_g(0,\mu,\nu)$
and $\sigma_{\mathrm{pr}}(\mu,\nu)=-\nu u_\sigma(0,\mu,\nu)$.
They cross transversely. The doubled singularity curve $\mu=\nu^2>0$
has quadratic contact with $\mu=0$. At every nonzero periodic root,
$g_s(s,\mu,\nu)=2\mu u_g(\mu,\mu,\nu)\ne0$, so the doubled cyclic
equations remain regular. Thus the punctured parabola records a singular
step on an existing doubled orbit, without a saddle-node of doubled orbits.

\begin{figure}[htbp]
\centering
\begin{minipage}{0.47\textwidth}
\centering
\begin{tikzpicture}[scale=1.5,>=stealth]
 \draw[->,gray] (-1.35,0)--(1.4,0) node[below,black] {$\nu$};
 \draw[->,gray] (0,-0.25)--(0,1.65) node[left,black] {$\mu$};
 \draw[thick,black!65] (-1.25,0)--(1.25,0);
 \draw[very thick,domain=-1.15:1.15,samples=80] plot (\x,{\x*\x});
 \fill (0,0) circle (1.3pt);
 \node at (0.85,1.5) {$\mu=\nu^2$};
 \node[below] at (-0.7,-0.06) {$\mu=0$};
\end{tikzpicture}
\par\small (a) Saddle-node interaction
\end{minipage}\hfill
\begin{minipage}{0.47\textwidth}
\centering
\begin{tikzpicture}[scale=1.5,>=stealth]
 \draw[->,gray] (-1.35,0)--(1.4,0) node[below,black] {$\nu$};
 \draw[->,gray] (0,-0.25)--(0,1.65) node[left,black] {$\mu$};
 \draw[thick,black!65] (-1.25,0)--(1.25,0);
 \draw[thick,black,dashed] (0,-0.2)--(0,1.45);
 \draw[very thick,domain=-1.15:1.15,samples=80] plot (\x,{\x*\x});
 \fill (0,0) circle (1.3pt);
 \node at (0.85,1.5) {$\mu=\nu^2$};
 \node[below] at (-0.7,-0.06) {$\mu=0$};
 \node[anchor=east] at (-0.06,0.8) {$\nu=0$};
\end{tikzpicture}
\par\small (b) Period-doubling interaction
\end{minipage}
\caption{The joint incidence diagrams in the coordinates of the two normal-form theorems. The horizontal line is the periodic bifurcation locus. In (a), the parabola consists of singular orbits on the two saddle-node branches. In (b), its punctured branches consist of singular doubled orbits; the dashed vertical line is the singularity locus on the primary branch. The axes are the local normal-form parameters.}
\label{fig:mixed-incidence}
\end{figure}
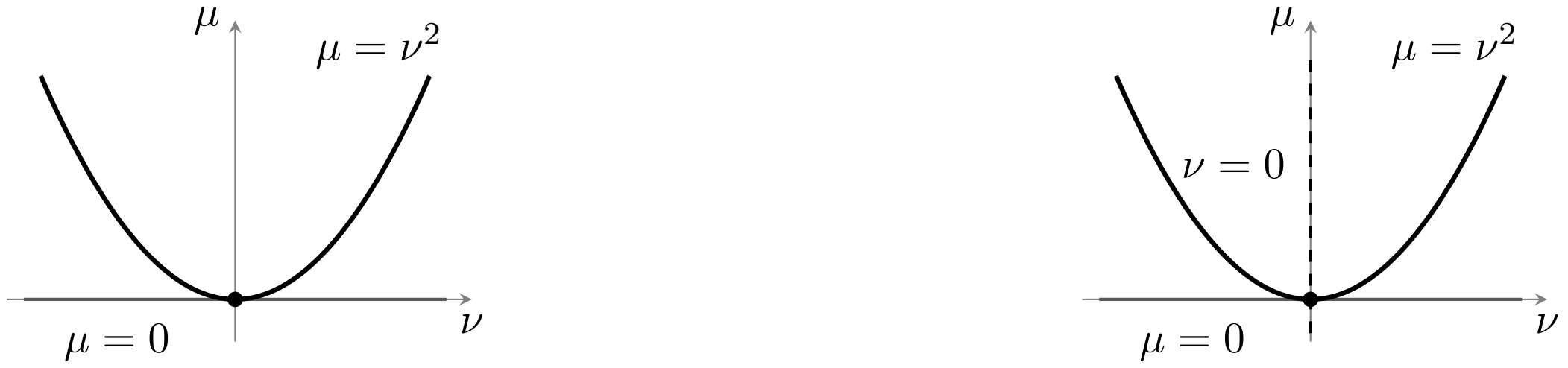

On the primary branch, $\sigma_+\sigma_-=\sigma_{\mathrm{pr}}^2$ has
zero differential at the collision. The separate signed markers are
therefore needed for the rank condition in \eqref{eq:mf-nondegeneracy}
and to identify which transition reaches the fold. The two mixed incidence diagrams are shown in Figure~\ref{fig:mixed-incidence}.

\subsection{Recognition coefficients}
The coefficient $a_3$ can be computed from the state derivatives of the block map $\mathcal D$.

In the following formula, $D^j\mathcal D$ denotes the $j$th derivative
with respect to the pair of state blocks $(R,S)$, with the parameter
fixed, evaluated at $(R,S,\lambda)=(X_*,X_*,0)$. Thus it is a
$j$-linear map from $(\R^{2d\ell})^j$ to $\R^{d\ell}$.
Define the vectors
\[
 \xi=(v_-,-v_-)\in\R^{2d\ell},\qquad
 q_2=D^2\mathcal D[\xi,\xi]\in\R^{d\ell},\qquad
 r_2=L_+^{-1}q_2\in\R^{d\ell}.
\]
Then
\begin{equation}\label{eq:mf-cubic}
 a_3=w_-^T\left\{
 \frac16D^3\mathcal D[\xi,\xi,\xi]
 -\frac12D^2\mathcal D[(r_2,r_2),\xi]\right\}.
\end{equation}
Differentiating \eqref{eq:mf-range} gives $v_s=z_s=0$, $v_{ss}=-r_2$, and $z_{ss}=0$ at the origin. Substitution in $g_{sss}$, followed by $w_-^TL_-=0$, gives \eqref{eq:mf-cubic}; this is the cubic calculation of the regular odd reduction \cite{Oliveira2026LS} with the same invertible range blocks.

We also need the variation of the linear coefficient
$\omega(0,\lambda)=g_s(0,\lambda)$. Let $L_-(\lambda)$ be the
operator $L_\ell(-1;\lambda)$ formed from the derivatives of the defining
relations along $X^{\mathrm{pr}}(\lambda)$. Thus its matrix entries
vary with both the orbit and the parameter. For a parameter direction
$\eta\in\R^2$,
\begin{equation}\label{eq:mf-linear-unfolding}
 D_\lambda\omega(0,0)\eta
   =w_-^T\left.\frac{\dd}{\dd t}\right|_{t=0}L_-(t\eta)v_-.
\end{equation}
Terms involving derivatives of the eliminated coordinates are
annihilated by $w_-^TL_-=0$. The orbit variation entering the derivative
of $L_-$ is
\[
 D_\lambda X^{\mathrm{pr}}(0)\eta
   =-L_+^{-1}D_\lambda\Phi(X_*,0)\eta.
\]
Thus \eqref{eq:mf-linear-unfolding} differentiates both partial
state derivatives of the relations and their evaluation points.

For the marker derivative, put $\widehat{v}_{-,k+1}=v_{-,k+1}$ if $k<\ell-1$ and $\widehat{v}_{-,k+1}=-v_{-,0}$ otherwise.

Differentiating at $(x_k^*,x_{k+1}^*;0)$ gives
\begin{equation}\label{eq:mf-marker-coupling}
 c_{\mathrm{PD}}=\operatorname{tr}\!\left(
 \adj(B_k)\bigl\{
 D_1(D_2G_k)[v_{-,k}]+D_2(D_2G_k)[\widehat{v}_{-,k+1}]
 \bigr\}\right).
\end{equation}
The vectors $v_{-,k}$ and $\widehat{v}_{-,k+1}$ belong to the cyclic anti-periodic mode and need not lie in $\ker B_k$.

\subsection{Invariance of the marked reduction}

A change of amplitude that respects block exchange has the form
$\widetilde s=\theta(s,\lambda)$, where $\theta$ is smooth,
$\theta(-s,\lambda)=-\theta(s,\lambda)$, and
$\theta_s(0,0)\ne0$. We call such a change equivariant.
Consider two equivariant Lyapunov--Schmidt reductions, written with
the same parameters and with the same distinguished transition.
Denote the second reduced equation and first marker by
$\widetilde g$ and $\widetilde\sigma_+$.
Let $\alpha$ and $\beta$ be smooth units, with $\alpha$ even in
$s$, and let $\gamma$ be a smooth scalar function. The equivalence
relevant to the periodic and singularity sets has the form
\begin{equation}\label{eq:mf-equivalence}
 \begin{aligned}
 \widetilde g(\theta(s,\lambda),\lambda)
   &=\alpha(s,\lambda)g(s,\lambda),\\
 \widetilde\sigma_+(\theta(s,\lambda),\lambda)
   &=\beta(s,\lambda)\sigma_+(s,\lambda)
       +\gamma(s,\lambda)g(s,\lambda).
 \end{aligned}
\end{equation}
The second marker obeys the corresponding
identity with $s$ replaced by $-s$. The term $\gamma g$ allows the
reconstructions to differ away from periodic solutions; it vanishes
on every solution under discussion.

\begin{proposition}\label{prop:mf-invariance}
The conditions \eqref{eq:mf-nondegeneracy} and the simultaneous
normal-form diagram are invariant under equivariant reduction
choices, regular defining equations, state coordinates, and starting
phase.
\end{proposition}

\begin{proof}
The reduction comparison in \cite{Oliveira2026LS} uses only invertibility of the range derivatives. To apply it here, use the left sides of the first range equations as coordinates $\mathfrak r\in\R^{2d\ell-1}$. Subtracting a smooth multiple of $\mathfrak r$ from the last equation puts the full system into the form $(\mathfrak r,g(s,\lambda))$. Write the left sides of the second range equations as
\[
 \mathcal J(s,\mathfrak r,\lambda)\mathfrak r
       +\mathfrak j(s,\mathfrak r,\lambda)g(s,\lambda),
\]
where $\mathcal J$ is a square matrix of size $2d\ell-1$ and
$\mathfrak j$ is a column vector of that length. At the reference
point $g=g_s=0$, so the derivative of this range system in
$(\mathfrak r,s)$ is $(\mathcal J,0)$. Its rank is $2d\ell-1$,
which makes $\mathcal J$ invertible there and hence nearby. On its zero set,
$\mathfrak r=T(s,\lambda)g(s,\lambda)$, with $T$ obtained by evaluating $-\mathcal J^{-1}\mathfrak j$ on the second range solution. Substitution into its last equation gives a unit times $g$, by the Schur complement of the invertible change of equations. The second amplitude defines $\theta$; exchange symmetry makes $\theta$ odd and the scalar unit $\alpha$ even.

The two reconstructed state vectors, when compared at these
corresponding amplitudes, differ by a smooth vector multiple of
$g$. Evaluating any smooth determinant marker on them therefore
changes its value by a smooth multiple of $g$. This proves
\eqref{eq:mf-equivalence} for the two reduction choices without
assuming that $g$ is a submersion.

For a regular change of equations $G_j\mapsto\Gamma_jG_j$, the
matrix $\Gamma_j$ is invertible and the marker on solutions is
multiplied by $\det\Gamma_j$. Away from solutions, differentiating
the product gives additional terms containing $G_j$. Along the
range reconstruction all components of the doubled map
are multiples of $g$, as shown above. Expansion of the determinant
therefore adds a multiple of $g$ and gives the same equivalence with
a unit $\beta$. State-coordinate changes likewise multiply the
determinants by units. They commute with block exchange when the
same charts are used at corresponding states. A change of starting
phase permutes the cyclic equations and transports the marked
transitions; it may exchange their plus and minus labels. In the
non-autonomous case the sequence of defining functions is shifted together with
the states.

To check the numerical conditions, write $\widetilde a_3$ and
$\widetilde c_{\mathrm{PD}}$ for the second reduction's cubic and amplitude
coefficients. Differentiating \eqref{eq:mf-equivalence} at the
reference point gives
\[
 \widetilde a_3\,\theta_s(0,0)^3=\alpha(0,0)a_3,
 \qquad
 \widetilde c_{\mathrm{PD}}\,\theta_s(0,0)=\beta(0,0)c_{\mathrm{PD}}.
\]
On $s=0$ the reduced linear coefficient is multiplied by
$\alpha(0,\lambda)/\theta_s(0,\lambda)$, and the marker by
$\beta(0,\lambda)$. Both original quantities vanish at the
reference parameter, so their parameter differentials are
multiplied by nonzero scalars. Their joint rank is preserved.
Finally, \eqref{eq:mf-equivalence} identifies the periodic and
singularity sets and preserves their parameter projections.
\end{proof}

\subsection{Relation to an ordinary flip}

For the spectral interpretation, impose \textup{(SP)}: every $A_j$ is invertible, and the primary backward return derivative $K(0)$ of \eqref{eq:HK} has simple eigenvalues $0,-1$, with all others nonzero and off the unit circle. Thus $d\geq2$. The degree-$d$ homogenisation of the cyclic determinant $\Delta_0$ has simple roots at infinity and at $-1$ by \eqref{eq:cyclic-backward-determinant}.

Assume \textup{(PD)} and \textup{(SP)}. On $\mu=0$, $\nu\ne0$, every primary transition is forward regular. The regular cyclic theory therefore applies: by \cite[Theorem~2.3]{Oliveira2026LS}, the condition $g_s=0$ identifies a multiplier $-1$, whose algebraic simplicity follows from \textup{(SP)} and continuity. The identities
\[
 g_s(0,\mu,\nu)=-\mu u_g(0,\mu,\nu),\qquad
 g_{sss}(0,0,\nu)=6u_g(0,0,\nu)\ne0
\]
give transverse crossing and cubic nondegeneracy by \cite[Proposition~6.4 and Theorem~6.5]{Oliveira2026LS}. Thus the punctured curve is an ordinary flip locus.

At its intersection with $\nu=0$ the forward derivative is not
defined. The mixed theorem describes how the closure of the regular
flip locus meets the fold of the evolution relation.

The doubled algebraic construction for explicit maps is established
in \cite{PeckhamKevrekidis1991}. It is used for singular discrete
relations in \cite[Section~5.2.1]{BW08}, and the regular cyclic
reduction is treated in \cite{Oliveira2026LS}. The two transition
markers retained here determine the additional joint incidence at an
arbitrary period. Formula~\eqref{eq:mf-linear-unfolding}
uses the derivative of the entire anti-periodic operator, including
both state derivatives and the continued orbit.

\section{Spectrum and stability of the doubled branch}
\label{sec:mixed-spectrum}

Assume \textup{(PD)} and \textup{(SP)}. We determine the escaping multiplier and the finite central multiplier on the primary and doubled branches.

\subsection{Determinant markers and the small backward eigenvalue}

Let $\lambda(\mu,\nu)$ be the inverse parameter change obtained in
Theorem~\ref{thm:mixed-flip}. Let $H_\lambda$ be the length-$\ell$
backward return of \eqref{eq:HK}, and define
\[
 K_{\rm pr}(\mu,\nu)
   =DH_{\lambda(\mu,\nu)}
        (x_0^{\mathrm{pr}}(\lambda(\mu,\nu))).
\]
On the doubled reconstruction define, for $0\leq j<2\ell$,
\[
 \begin{aligned}
 A_j^{(2)}(s,\mu,\nu)
  &=D_1G_{j\bmod\ell}
      (Y_j(s,\mu,\nu),Y_{j+1}(s,\mu,\nu);\lambda(\mu,\nu)),\\
 B_j^{(2)}(s,\mu,\nu)
  &=D_2G_{j\bmod\ell}
      (Y_j(s,\mu,\nu),Y_{j+1}(s,\mu,\nu);\lambda(\mu,\nu)),\\
 C_j^{(2)}(s,\mu,\nu)
  &=-\bigl(A_j^{(2)}(s,\mu,\nu)\bigr)^{-1}B_j^{(2)}(s,\mu,\nu),\\
 K_2(s,\mu,\nu)&=C_0^{(2)}(s,\mu,\nu)\cdots
                       C_{2\ell-1}^{(2)}(s,\mu,\nu).
 \end{aligned}
\]
The cyclic convention is $Y_{2\ell}=Y_0$. All inverses in this
definition exist by \textup{(SP)} and continuity.
The matrix $K_2$ is a smooth algebraic product even where $g\ne0$.
At a doubled periodic solution, where $s^2=\mu>0$, it is the
derivative of the actual length-$2\ell$ backward return.
At $s=0$ the two blocks coincide, so
$K_2(0,\mu,\nu)=K_{\rm pr}(\mu,\nu)^2$.

Let $\kappa_{\rm pr}(\mu,\nu)$ and $\kappa_2(s,\mu,\nu)$ continue the simple zero eigenvalues of $K_{\rm pr}(0,0)$ and $K_2(0,0,0)=K_{\rm pr}(0,0)^2$. On the regular branches set
\[
 \rho_{\rm pr}=\kappa_{\rm pr}^{-1},\qquad
 \rho_2=\bigl(\kappa_2|_{s^2=\mu}\bigr)^{-1}.
\]
Restrictions to $s^2=\mu$ below apply only to even functions of $s$. They are independent of the sign of $s$, smooth for $\mu\ge0$, and extended smoothly across $\mu=0$ by the even-function factorisation. The actual doubled branch has $\mu>0$. The coefficients $e_{\rm pr},c_\infty$ below are functions of $(\mu,\nu)$ with this extension convention.

\begin{proposition}\label{prop:mixed-poles}
Under \textup{(PD)} and \textup{(SP)}, a smooth unit $e_{\rm pr}$
and a positive smooth function $c_\infty$ satisfy
\begin{align}
 \kappa_{\rm pr}(\mu,\nu)&=e_{\rm pr}(\mu,\nu)\nu,
 &\rho_{\rm pr}&=\frac{e_{\rm pr}^{-1}}{\nu},
 \label{eq:primary-flip-pole-law}\\
 \kappa_2\big|_{s^2=\mu}&=c_\infty(\mu,\nu)(\nu^2-\mu),
 &\rho_2&=\frac{c_\infty(\mu,\nu)^{-1}}{\nu^2-\mu}.
 \label{eq:doubled-pole-law}
\end{align}
The reciprocal formulas apply to the forward-regular orbits.
\end{proposition}
\begin{proof}
The sole possibly vanishing determinant on the primary branch is
$\sigma_{\mathrm{pr}}(\mu,\nu)=-\nu u_\sigma(0,\mu,\nu)$. Applying
Lemma~\ref{lem:marker-spectral-product} to the primary product gives
$\kappa_{\rm pr}=e_{\rm pr}\nu$, with $e_{\rm pr}$ a unit.
The possibly vanishing determinants in the doubled product are
$\sigma_+$ and $\sigma_-$. The same factorisation and
\eqref{eq:mf-normal-plus}--\eqref{eq:mf-normal-minus} therefore give
\[
 \kappa_2(s,\mu,\nu)=e_2(s,\mu,\nu)(\nu^2-s^2)
\]
for a smooth unit $e_2$.

To check that this coefficient descends to the unlabelled doubled
orbit, define the half-products
\[
 P_2=C_0^{(2)}\cdots C_{\ell-1}^{(2)},\qquad
 Q_2=C_\ell^{(2)}\cdots C_{2\ell-1}^{(2)}.
\]
Block exchange replaces $P_2Q_2$ by $Q_2P_2$. These two square
matrices have the same characteristic polynomial, including when
one is singular. Since $\kappa_2$ is the unique eigenvalue near
zero, it is even in $s$. The identity above makes $e_2$ even wherever
$\nu^2\ne s^2$, and hence everywhere by continuity. Write
$e_2(s,\mu,\nu)=\widehat e_2(s^2,\mu,\nu)$ with
$\widehat e_2$ smooth, and set
$c_\infty(\mu,\nu)=\widehat e_2(\mu,\mu,\nu)$.
This gives the doubled factorisation.

At $s=0$ the small eigenvalue of $K_2=K_{\rm pr}^2$ is
$\kappa_{\rm pr}^2$. For $\nu\ne0$ it follows that
$e_2(0,\mu,\nu)=e_{\rm pr}(\mu,\nu)^2$.
By continuity the same equality holds at $\nu=0$. In particular,
$c_\infty(0,0)>0$, so it remains positive after the neighbourhood
has been reduced. Taking reciprocals proves the asserted pole laws
on the regular branches.
\end{proof}

\subsection{The finite central multiplier}

Let $\eta_{\rm pr}(\mu,\nu)$ be the backward eigenvalue of
$K_{\rm pr}$ near $-1$, and let $\eta_2(s,\mu,\nu)$ be the
backward eigenvalue of $K_2$ near $1$. They are simple and smooth:
under \textup{(SP)}, squaring $K_{\rm pr}(0,0)$ gives a simple
eigenvalue $1$, because no other eigenvalue equals $1$ or $-1$.
The block-exchange argument used above makes $\eta_2$ even in $s$.
Define
\[
 m_{\rm pr}=\eta_{\rm pr}^{-1},\qquad
 m_2=\bigl(\eta_2|_{s^2=\mu}\bigr)^{-1}.
\]
These spectral functions remain smooth at singular evolution, but represent forward stability only on regular orbits. Write $d_{\rm pr},d_2$ for the coefficient functions of $(\mu,\nu)$ in the following identities, using the same extension convention.

\begin{proposition}\label{prop:flip-central-multipliers}
Under \textup{(PD)} and \textup{(SP)}, smooth units $d_{\rm pr}$
and $d_2$ satisfy
\begin{equation}\label{eq:flip-central-laws}
 m_{\rm pr}+1=d_{\rm pr}(\mu,\nu)\mu,
 \qquad m_2-1=d_2(\mu,\nu)\mu,
 \qquad d_2(0,\nu)=4d_{\rm pr}(0,\nu).
\end{equation}
\end{proposition}
\begin{proof}
Recall that $\Psi$ is the doubled map and $Y$ is its
range reconstruction. At fixed parameter, Schur elimination of the
invertible derivative of the range equations gives a smooth unit
$a_{\det}$ such that
\[
 \det D_{(R,S)}\Psi(Y(s,\mu,\nu),\lambda(\mu,\nu))
       =a_{\det}(s,\mu,\nu)g_s(s,\mu,\nu).
\]
The identity holds on the whole range reconstruction, including $g\ne0$, because the eliminated range derivative remains invertible.

The cyclic determinant identity
\eqref{eq:cyclic-backward-determinant}, now for length $2\ell$ and
twist $1$, expresses the same determinant as
$\bigl(\prod_{j=0}^{2\ell-1}\det A_j^{(2)}\bigr)\det(I-K_2)$.
All factors of $\det(I-K_2)$ except $1-\eta_2$ are nonzero nearby.
Consequently there is a smooth unit $b_c$ such that
\begin{equation}\label{eq:central-schur-identity}
 1-\eta_2(s,\mu,\nu)=b_c(s,\mu,\nu)g_s(s,\mu,\nu).
\end{equation}
Both $\eta_2$ and $g_s$ are even, and $\partial_\mu g_s(0,0,0)=-u_g(0,0,0)\ne0$. Thus $g_s\ne0$ on a dense subset nearby. Division there, followed by continuity, makes $b_c$ even.

Set $b_{c,0}(\mu,\nu)=b_c(0,\mu,\nu)$ and
$u_{g,0}(\mu,\nu)=u_g(0,\mu,\nu)$. At $s=0$,
$\eta_2=\eta_{\rm pr}^2$ and $g_s=-\mu u_{g,0}$. Factoring
$1-\eta_{\rm pr}^2$ in \eqref{eq:central-schur-identity} and
then dividing by the nonzero $\eta_{\rm pr}$ gives
\[
 m_{\rm pr}+1=d_{\rm pr}\mu,
 \qquad
 d_{\rm pr}=-\frac{b_{c,0}u_{g,0}}
                   {\eta_{\rm pr}(1-\eta_{\rm pr})}.
\]
At a nonzero periodic root, $g_s=2\mu u_g(\mu,\mu,\nu)$.
Dividing \eqref{eq:central-schur-identity} by $\eta_2$ therefore
gives
\[
 m_2-1=d_2\mu,
 \qquad
 d_2=\left.
      \frac{2b_c(s,\mu,\nu)u_g(\mu,\mu,\nu)}
           {\eta_2(s,\mu,\nu)}\right|_{s^2=\mu}.
\]
Evenness makes this restriction smooth up to $\mu=0$, with a smooth
extension across that boundary. At $\mu=0$,
\eqref{eq:central-schur-identity} implies
$\eta_{\rm pr}=-1$ and $\eta_2=1$. The two formulas then give
$d_{\rm pr}=b_{c,0}u_{g,0}/2$ and $d_2=2b_{c,0}u_{g,0}$, proving the factor $4$.
Every factor in these expressions is nonzero nearby.
\end{proof}

The two units $d_{\rm pr}$ and $d_2$ have the same sign near the
origin. Thus, when $\mu>0$ and both orbits are forward regular,
exactly one of the primary and doubled orbits is attracting in its
finite central direction. Near $-1$, attraction means
$m_{\rm pr}+1>0$; near $1$, it means $m_2-1<0$.
This is a statement about one direction. Both full orbits have an
unstable direction from their small backward eigenvalue, since its
nonzero reciprocal has modulus greater than one near the collision.

The doubled pole law also fixes its sign:
$\rho_2>0$ for $\mu<\nu^2$ and $\rho_2<0$ for $\mu>\nu^2$
on the regular doubled branch. By contrast, traversing the primary
orbit twice squares its large multiplier:
\[
 \rho_{\rm pr}^2=e_{\rm pr}^{-2}\nu^{-2}.
\]
The repeated primary sequence has two copies of the same singular
transition. They produce one multiplier with a pole of order two,
not two independent multipliers at infinity.

Finally, a simple zero eigenvalue of the repeated backward derivative
does not by itself make the repeated return a fold. To see the
obstruction at the collision, let $H=H_0$ be the length-$\ell$
backward return at the reference parameter and put $K=DH(x_0^*)$.
Choose nonzero right and left zero-eigenvectors $v_K,w_K\in\R^d$,
so $Kv_K=0$ and $w_K^TK=0$. Algebraic simplicity of zero implies
$\ker K^2=\ker K$ and $\ker(K^2)^T=\ker K^T$.
Since $H(x_0^*)=x_0^*$, the chain rule gives
\[
 \begin{aligned}
 w_K^TD^2(H\circ H)(x_0^*)[v_K,v_K]
 &=w_K^TD^2H(x_0^*)[Kv_K,Kv_K]
       +w_K^TKD^2H(x_0^*)[v_K,v_K]\\
 &=0.
 \end{aligned}
\]
The quadratic coefficient required for a fold therefore vanishes
for the repeated return at the collision. On the punctured parabola
$\mu=\nu^2>0$, however, the doubled orbit has exactly one singular
transition. The single-transition return-fold result of
Proposition~\ref{prop:returnfold} applies there.

\section{Examples and realisation at prescribed periods}
\label{sec:examples}

The return-model examples are written in the form
$F(x,y;\lambda)=h_\lambda(y)-x$, with the appropriate phase index in
the non-autonomous case. Thus $x=h_\lambda(y)$ determines the
predecessor of $y$, and $h_\lambda$ is a backward map. The two state
derivatives are $A=-I_d$ and $B=Dh_\lambda(y)$, where $I_d$ is the
identity matrix in dimension $d$. In particular, the relation is
regular even when forward evolution fails.

\subsection{A scalar two-cycle and a non-autonomous four-cycle}
\label{subsec:autonomous-two-example}

Let the state variables $x,y$ and the parameter $\lambda$ be real.
Define
\begin{equation}\label{eq:autonomous-two-relation}
 h_\lambda(y)=1-y^2+\lambda y(y-1),
 \qquad F(x,y;\lambda)=h_\lambda(y)-x.
\end{equation}
Since $h_\lambda(1)=0$ and $h_\lambda(0)=1$, the cyclic state
sequence $(x_0,x_1)=(0,1)$ is an orbit of period two for every
$\lambda$. Here $A_0=A_1=-1$. Evaluating $h_\lambda'$ at the
successor of each transition gives the next-state derivatives
$B_0(\lambda)=h_\lambda'(1)$ and $B_1(\lambda)=h_\lambda'(0)$. The cyclic derivative
$L(\lambda)$ and the derivative $K(\lambda)$ of the two-step backward return
$h_\lambda\circ h_\lambda$, evaluated at $x_0=0$, are therefore
\begin{equation}\label{eq:autonomous-two-matrices}
 B_0(\lambda)=\lambda-2,\qquad B_1(\lambda)=-\lambda,\qquad
 L(\lambda)=\begin{pmatrix}-1&\lambda-2\\-\lambda&-1\end{pmatrix},
 \qquad K(\lambda)=\lambda(2-\lambda).
\end{equation}
The determinant $\det L(\lambda)=(1-\lambda)^2$ is nonzero near $\lambda=0$.
Consequently the displayed orbit is the unique nearby solution of
the equations for closure after two steps. Nevertheless its transition from $1$ to
$0$ is forward singular at $\lambda=0$. The corresponding marker is
$-\lambda$, which crosses zero with nonzero derivative, and
$\partial_y^2F(1,0;0)=-2$. This is an evolution fold on a smoothly
continued periodic orbit.

For $\lambda$ near zero, let $y_c(\lambda)$ be the critical point of
$h_\lambda$ and let $x_c(\lambda)=h_\lambda(y_c(\lambda))$ be its
critical value. Direct calculation gives
\[
 y_c(\lambda)=-\frac{\lambda}{2(1-\lambda)},\qquad
 x_c(\lambda)=1+\frac{\lambda^2}{4(1-\lambda)}.
\]
Completing the square rewrites $F=0$ as
\begin{equation}\label{eq:autonomous-two-square}
 x=x_c(\lambda)-(1-\lambda)(y-y_c(\lambda))^2.
\end{equation}
For a present state $x$ near $1$, there are two nearby successors
when $x<x_c(\lambda)$, one when $x=x_c(\lambda)$, and none when
$x>x_c(\lambda)$. At the periodic present state $x=1$, the two
successors are $0$ and $-\lambda/(1-\lambda)$. Their merger at
$\lambda=0$ does not create another nearby two-cycle.

Here (S) holds since $K(0)=0$. Proposition~\ref{prop:uniform-all-period-isolation} excludes all additional multiple-period orbits in a fixed local itinerary. For nonzero $\lambda$, the forward multiplier is
\begin{equation}\label{eq:autonomous-two-forward}
 M(\lambda)=K(\lambda)^{-1}=\frac{1}{\lambda(2-\lambda)}.
\end{equation}

Now prescribe period $p=2$ for the defining sequence. Use the same polynomial
$h_\lambda$ and define the two phase-dependent backward maps and
relations by
\begin{equation}\label{eq:forced-two-relations}
 h_{0,\lambda}(y)=h_\lambda(y-2),\qquad
 h_{1,\lambda}(y)=y+2,\qquad
 F_j(x,y;\lambda)=h_{j,\lambda}(y)-x\quad(j=0,1).
\end{equation}
The equation at time $n$ uses $F_{n\bmod2}$. One backward map is
quadratic and the other is affine, so the defining sequence has
period two near $\lambda=0$. The cyclic sequence $(0,3,1,2)$ is an
orbit of period four. Its successive next-state derivatives are
\[
 (B_0,B_1,B_2,B_3)=(\lambda-2,1,-\lambda,1).
\]
Exactly its transition from $1$ to $2$ is singular at $\lambda=0$.
The two-step backward return from phase zero to phase zero is
$h_{0,\lambda}\circ h_{1,\lambda}=h_\lambda$. Consequently the
four-step backward return is $h_\lambda\circ h_\lambda$; its
derivative, cyclic determinant, and forward multiplier are
$\lambda(2-\lambda)$, $(1-\lambda)^2$, and
$[\lambda(2-\lambda)]^{-1}$, respectively. Proposition~\ref{prop:uniform-all-period-isolation} applies to the four-step return. Here $N=2$, so the
orbit period $Np=4$ is strictly larger than the period of the defining sequence.

\subsection{A fold meeting a saddle-node}
\label{subsec:saddle-example}

Let $z=(u,v)\in\R^2$ be a state and let $(\mu,\nu)\in\R^2$ be the
parameter pair. Define the backward map $H_{\mu,\nu}$ and its
implicit relation by
\begin{equation}\label{eq:saddle-explicit-model}
 H_{\mu,\nu}(u,v)
   =\bigl(u+u^2-\mu,\ v^2+(u-\nu)v\bigr),\qquad
 F(x,y;\mu,\nu)=H_{\mu,\nu}(y)-x,
 \quad x,y\in\R^2.
\end{equation}
The second fixed-point equation is $v(v+u-\nu-1)=0$. The second factor is nonzero in
a sufficiently small fixed neighbourhood of the origin in state and
parameter space. The fixed-point equations there are therefore
exactly
\begin{equation}\label{eq:saddle-example-solutions}
 u^2=\mu,\qquad v=0.
\end{equation}
At such a fixed point, write $K=DH_{\mu,\nu}(u,0)$ for the backward
return derivative. Since a fixed point has a single transition,
$B=K$, and
\begin{equation}\label{eq:saddle-example-spectrum}
 K=DH_{\mu,\nu}(u,0)
   =\operatorname{diag}(1+2u,u-\nu),\qquad
 \det B=(1+2u)(u-\nu).
\end{equation}
Removing the unit $1+2u$ gives the marker $u-\nu$, so the singular periodic set projects onto $\mu=\nu^2$. At $u=\nu$, the kernel of $B$ is the $v$-axis and the second component has quadratic derivative $2$ in that direction. The mixed rank condition is the invertibility of $D_{(u,\mu,\nu)}(u^2-\mu,2u,u-\nu)$ at zero. Thus (SN) holds. The backward eigenvalues at the collision are $1$ and $0$, and the regular forward multipliers are
\begin{equation}\label{eq:saddle-example-forward}
 \frac{1}{1+2u},\qquad \frac{1}{u-\nu}.
\end{equation}
On $u=\sqrt\mu$ the finite mode is attracting, and on $u=-\sqrt\mu$ it is repelling. The pole mode makes both regular orbits unstable near the collision, as in Proposition~\ref{prop:ms-spectrum}.

\subsection{A mixed interaction without backward regularity}
\label{subsec:singular-backward-example}
Consider the scalar $2$-periodic relations
\begin{equation}\label{eq:singular-backward-relations}
 F_0(x,y;\mu,\nu)=(y-1)^2-x,\qquad
 F_1(x,y;\mu,\nu)=2(x-1)^2-y+\nu(x-1)-\mu.
\end{equation}
At $(\mu,\nu)=(0,0)$, the state sequence $(0,1)$ has period two, and
\[
 (A_0,B_0,A_1,B_1)=(-1,0,0,-1),\qquad
 L=\begin{pmatrix}-1&0\\-1&0\end{pmatrix}.
\]
Both transitions are relation regular. The first is an evolution fold, with $(F_0)_{yy}=2$; the second has no smooth local predecessor at the collision, since $A_1=0$ and $(F_1)_{xx}=4$.

Choose cyclic kernel vectors $v_L=(0,1)^T$ and $w_L=(-1,1)^T/\sqrt2$, and write $t=x_1-1$. The range equation gives $x_0=(3t^2+\nu t-\mu)/2$. The reduced equation and the marker of the first transition are
\begin{equation}\label{eq:singular-backward-reduction}
 h(t,\mu,\nu)=\frac{t^2+\nu t-\mu}{\sqrt2},\qquad
 \sigma(t,\mu,\nu)=2t.
\end{equation}
At zero, $h_{tt}=\sqrt2$, $\sigma_t=2$, and
$\det D_{(t,\mu,\nu)}(h,h_t,\sigma)=-1$, so \textup{(SN)} holds. The periodic solutions satisfy $x_0=t^2$ and $t^2+\nu t-\mu=0$. Their fold locus is $\mu=-\nu^2/4$, while their singular-transition locus is $\mu=0$. Equivalently, $u=t+\nu/2$, $\widehat\mu=\mu+\nu^2/4$ and $\widehat\nu=\nu/2$ give
$h=(u^2-\widehat\mu)/\sqrt2$ and $\sigma=2(u-\widehat\nu)$.
Thus Theorem~\ref{thm:ms-normal} describes the interaction even though the second transition cannot be solved backward through the collision.

\subsection{A fold meeting period doubling}
\label{subsec:flip-example}

Let $(u,v)$ be the state coordinates and $(\mu,\nu)$ the parameters.
Define a scalar map $\varphi_\mu$, a planar backward map
$H_{\mu,\nu}$, and its implicit relation by
\begin{equation}\label{eq:flip-explicit-model}
 \begin{aligned}
 \varphi_\mu(u)&=-(1+\mu)u+u^3,\\
 H_{\mu,\nu}(u,v)&=\bigl(\varphi_\mu(u),v^2+(\nu+u)v\bigr),\\
 F(x,y;\mu,\nu)&=H_{\mu,\nu}(y)-x,
 \qquad x,y\in\R^2.
 \end{aligned}
\end{equation}
The fixed-point equations are $u(u^2-2-\mu)=0$ and
$v(v+\nu+u-1)=0$. Thus the origin is the only sufficiently small
fixed point. Let $K$ be its backward derivative and, when it exists,
let $M$ be its forward derivative. They are
\begin{equation}\label{eq:flip-example-primary-spectrum}
 K=\operatorname{diag}(-1-\mu,\nu),\qquad
 M=\operatorname{diag}\bigl(-(1+\mu)^{-1},\nu^{-1}\bigr)
 \quad(\nu\ne0).
\end{equation}
The primary orbit has a singular transition on $\nu=0$ and a finite
multiplier $-1$ on $\mu=0$. These two conditions have independent
parameter derivatives. At their intersection the kernel of $B$ is
the $v$-axis, and the second derivative of the second component in
this direction is $2$. The evolution singularity is a fold.

The recognition coefficients can be checked directly. With
$v_-=w_-=(1,0)^T$, the doubled range reconstruction is
$R=(s,0)$, $S=(-s,0)$, and
\[
 g(s,\mu,\nu)=s(\mu-s^2),\qquad
 \sigma_+(s,\mu,\nu)=(-1-\mu+3s^2)(\nu-s).
\]
At the collision, $a_3=-1$, $c_{\mathrm{PD}}=1$, and
$\det D_{(\mu,\nu)}\bigl(\mu,-(1+\mu)\nu\bigr)=-1$.
Thus \textup{(PD)} holds; the simple backward eigenvalues $-1,0$
also give \textup{(SP)}.

For $\mu>0$, set $u_*=\sqrt\mu$. The states $(u_*,0),(-u_*,0)$ form the unique small two-cycle. Indeed, choose $\theta\in(0,1)$ bounding $|v+\nu+u|$ locally. Iteration of the second coordinate around any periodic orbit forces $v=0$. If $u,w$ are distinct first coordinates of a two-cycle, adding $w=\varphi_\mu(u)$ and $u=\varphi_\mu(w)$ gives
\[
 (u+w)(2+\mu-u^2+uw-w^2)=0.
\]
The second factor is nonzero near zero, hence $w=-u$ and $u^2=\mu$.

The two transition determinants on this orbit are
\begin{equation}\label{eq:flip-example-transition-markers}
 \det DH_{\mu,\nu}(\pm u_*,0)=(-1+2\mu)(\nu\pm u_*).
\end{equation}
The factor $-1+2\mu$ is a unit near $\mu=0$. Thus the two markers
are $\nu+u_*$ and $\nu-u_*$. The doubled orbit has a singular
transition precisely on $\mu=\nu^2>0$. Exactly one marker vanishes
there, so exactly one of the two transitions is singular. Its
second derivative in the kernel direction is again $2$.

Let $K_2$ denote the derivative of the two-step backward return at
either state of the two-cycle. Let $m_{\mathrm c}$ and
$\rho_\infty$ denote its two forward multipliers where both
transitions are regular; the subscripts distinguish the multiplier
near $1$ from the multiplier that diverges. Multiplication of the
two diagonal step derivatives gives
\begin{align}
 K_2&=\operatorname{diag}\bigl((-1+2\mu)^2,\nu^2-\mu\bigr),
 \label{eq:flip-example-doubled-spectrum}\\
 m_{\mathrm c}&=\frac{1}{(-1+2\mu)^2},\qquad
 \rho_\infty=\frac{1}{\nu^2-\mu}.
 \label{eq:flip-example-doubled-forward}
\end{align}
The escaping denominator is the marker product. On the repeated primary orbit the multiplier is instead $\nu^{-2}$, in agreement with the order-two pole of Section~\ref{sec:mixed-spectrum}. On the doubled branch, $m_{\mathrm c}>1$ and $|\rho_\infty|>1$ near the collision, so both directions are unstable.

\subsection{Autonomous and non-autonomous realisation of the return models}
\label{subsec:prescribed-period-realisation}

Let $H_\lambda$ be a smooth family defined near $0\in\R^d$, with $H_0(0)=0$ and $\lambda$ near $0\in\R^m$. We realise it as a backward return at any prescribed period.

Fix an integer $q\geq1$. Let $e_1=(1,0,\ldots,0)\in\R^d$ and
choose the centres $c_j=4je_1$ for $0\leq j<q$. Centre indices are
read modulo $q$. Choose $r\in(0,1)$ small enough that
$H_\lambda$ is defined for $\|z\|<r$ and small $\lambda$.
The balls of radius $r$ about the centres are pairwise disjoint.
On smaller balls define a single backward map $h_\lambda$ by
\begin{equation}\label{eq:autonomous-chart-realisation}
 h_\lambda(c_1+z)=c_0+H_\lambda(z),\qquad
 h_\lambda(c_{j+1}+z)=c_j+z\quad(1\leq j<q).
\end{equation}
For $q=1$, only the first formula is used and $c_1=c_0$.
Smooth cutoff functions supported in the disjoint balls extend these
prescriptions to a smooth map on $\R^d$. Since $H_0(0)=0$, further
restriction of the state and parameter neighbourhoods ensures that
all compositions below stay in the prescribed charts.

Starting at $c_0+z$, the backward itinerary traverses $q-1$ translation charts before applying $H_\lambda$. Its return is exactly $c_0+z\mapsto c_0+H_\lambda(z)$. Set $F(x,y;\lambda)=h_\lambda(y)-x$.

\begin{proposition}\label{prop:prescribed-period-realisation}
The reference orbit has period $q$. Nearby fixed points
and two-cycles of $H_\lambda$ correspond to $q$- and
$2q$-orbits of the relation, respectively.
\end{proposition}
\begin{proof}
The centres give the reference orbit. A fixed point $z$ of $H_\lambda$ gives $(c_0+z,\ldots,c_{q-1}+z)$, of period $q$ by disjointness of the charts. A two-cycle $z,w$ alternates between these coordinates on successive visits to $c_0$. Its period is a multiple of $q$ and differs from $q$, hence is $2q$. Conversely, sampling a periodic sequence at that chart gives the corresponding orbit of $H_\lambda$.
\end{proof}

Every present-state derivative is $-I_d$, and every added next-state derivative is $I_d$. The construction therefore preserves the return spectrum, the singular determinant and its fold coefficient. It applies to both planar models. The scalar model $H_\lambda(z)=z^2+\lambda z$ likewise gives a persistent $q$-orbit with marker $\lambda$ and forward multiplier $\lambda^{-1}$.

For the non-autonomous construction, fix a period $p\geq2$
for the defining sequence and an integer $N\geq1$. Apply the preceding construction
with $q=N$, so its centres are $c_k=4ke_1$ for $0\leq k<N$.
Set $b=4Ne_1$. The translated centres
\[
 c_k+jb=4(k+jN)e_1,
 \qquad 0\leq k<N,\quad 0\leq j<p,
\]
are pairwise distinct, with mutual distances at least $4$.
Their balls of radius $r<1$ are therefore disjoint. Define the
phase-dependent backward maps and relations by
\begin{equation}\label{eq:forced-general-relations}
 \begin{aligned}
 h_{0,\lambda}(y)&=h_\lambda(y-(p-1)b),\\
 h_{j,\lambda}(y)&=y+b\quad(1\leq j<p),\\
 F_j(x,y;\lambda)&=h_{j,\lambda}(y)-x\quad(0\leq j<p).
 \end{aligned}
\end{equation}
The equation at time $n$ uses $F_{n\bmod p}$. Composition in backward
order gives
$h_{0,\lambda}\circ h_{1,\lambda}\circ\cdots\circ h_{p-1,\lambda}
=h_\lambda$. Thus the $p$-step backward return at phase zero is the
autonomous map just constructed. In each example used here,
$h_{0,\lambda}$ is nonlinear and the remaining phase maps are affine.
Exactly one function in each period of the defining sequence is
nonlinear, so the period of this sequence is $p$.

Let $(z_0,\ldots,z_{N-1})$ be an $N$-orbit of the
autonomous relation, so $h_\lambda(z_{k+1})=z_k$ with indices
modulo $N$. Define its lift by
\begin{equation}\label{eq:forced-general-orbit}
 x_{kp}=z_k,\qquad
 x_{kp+j}=z_{k+1}+(p-j)b
 \quad(0\leq k<N,\ 1\leq j<p).
\end{equation}
At phase zero,
$h_{0,\lambda}(x_{kp+1})=h_\lambda(z_{k+1})=z_k$;
the remaining transitions follow by subtracting $b$.
The lifted orbit has period $Np$: phase recurrence requires a multiple of $p$ steps, and the phase-zero subsequence has period $N$. The same construction lifts $2N$-orbits to $2Np$-orbits.

All added transitions have identity next-state derivative. The
$Np$-step backward return is $H_\lambda$, and the $2Np$-step return
on a doubled orbit is $H_\lambda\circ H_\lambda$. The singularity
curves and multiplier formulas of the original return models are
therefore unchanged.

Finally, the planar mixed examples extend to any state dimension
$d>2$. Let $w\in\R^{d-2}$ be the additional state coordinate and
choose an invertible real $(d-2)\times(d-2)$ matrix $C_\perp$ with
no eigenvalue of modulus one. Replace the planar backward map by
\[
 \widehat H_\lambda(u,v,w)
   =\bigl(H_\lambda(u,v),C_\perp w\bigr).
\]
The additional coordinates of a nearby periodic orbit must be zero,
since no positive power of $C_\perp$ has eigenvalue $1$.
The two interacting directions are unchanged. The other forward
multipliers are the eigenvalues of $C_\perp^{-1}$ on the primary
orbit and of $C_\perp^{-2}$ on a doubled orbit.

\subsection{Applied implicit relations}\label{subsec:applied-relations}
The economic and machining calculations concern the single evolution fold. The backward-Euler examples also include a two-parameter interaction between an evolution fold and period doubling.

\paragraph{Labour-market equilibrium.}
In the deterministic search-and-matching model of Lazaryan and Lubik \cite{LazaryanLubik2019} with risk-neutral households, market tightness is the ratio of vacancies to unemployment. Write $z$ for its square root. Let $\beta_{\rm L},\rho_{\rm L},\eta_{\rm L}\in(0,1)$ denote the discount factor, separation fraction and worker bargaining weight, respectively. Let $m_{\rm L},s_{\rm L},\kappa_{\rm L}>0$ denote matching efficiency, productivity net of unemployment income and the vacancy-posting cost. With matching elasticity $1/2$, define
\[
 a_{\rm L}=\beta_{\rm L}(1-\rho_{\rm L}),\qquad
 b_{\rm L}=\eta_{\rm L}m_{\rm L},\qquad
 C_{\rm L}=\frac{(1-\eta_{\rm L})m_{\rm L}s_{\rm L}}{\kappa_{\rm L}}.
\]
The equilibrium relation for successive values $x,y>0$ of $z$ is
\begin{equation}\label{eq:application-labour}
 F_{\rm L}(x,y;C_{\rm L})
 =x-a_{\rm L}(C_{\rm L}+y-b_{\rm L}y^2)=0.
\end{equation}
The following parameter choice and singularity calculation are ours. Define
\[
 z_*=\frac{1}{2b_{\rm L}},\qquad
 C_*=\frac{1}{2a_{\rm L}b_{\rm L}}-\frac{1}{4b_{\rm L}},
 \qquad \lambda=C_{\rm L}-C_*.
\]
At $(x,y,C_{\rm L})=(z_*,z_*,C_*)$,
\[
 F_{\rm L}=0,\qquad A=1,\qquad B=0,\qquad
 (F_{\rm L})_{yy}=2a_{\rm L}b_{\rm L}>0,\qquad L=1.
\]
The steady equilibrium continues as $z=z(C_{\rm L})$, with $z'(C_*)=a_{\rm L}$. Its marker is
$\sigma(C_{\rm L})=a_{\rm L}(2b_{\rm L}z(C_{\rm L})-1)$, so
\begin{equation}\label{eq:application-labour-crossing}
 \sigma'(C_*)=2a_{\rm L}^2b_{\rm L}>0.
\end{equation}
All conditions (E) hold, and the backward derivative is zero. Theorems~\ref{thm:contact} and~\ref{thm:simple-infinite-multiplier} therefore give quadratic contact of current tightness with the successor boundary and the forward multiplier
$-[2a_{\rm L}^2b_{\rm L}\lambda]^{-1}+O(1)$.

This event lies in the admissible probability region. Take $m_{\rm L}=1/2$, $\eta_{\rm L}=3/4$, $\beta_{\rm L}=99/100$, $\rho_{\rm L}=1/10$, $s_{\rm L}=3/5$ and $\kappa_{\rm L}=8019/88720$. Then $a_{\rm L}=891/1000$, $b_{\rm L}=3/8$, $z_*=4/3$ and $C_*=2218/2673$. At tightness $z_*^2=16/9$, the job-finding and vacancy-filling probabilities are $2/3$ and $3/8$, respectively. Thus the equilibrium relation loses a smooth local rule for future tightness while the stationary equilibrium persists. Nearby current states admit two continuations or none.

\paragraph{Machining.}
Let $x,y$ be successive surface errors, let $a_0$ be prescribed radial immersion, and let $f(a)$ be the smooth surface displacement response at actual immersion $a$, in the interior of the cutting regime. The milling relation of Kiss and Bachrathy \cite{KissBachrathy2022} is
\begin{equation}\label{eq:application-milling}
 F_{\rm M}(x,y;a_0)=y-f(a_0+x-y)=0.
\end{equation}
Its fixed point is $x=y=f(a_0)$. Along it,
\[
 A=-f'(a_0),\qquad B=1+f'(a_0),\qquad L=1.
\]
At an immersion $a_*$ with $f'(a_*)=-1$ and $f''(a_*)\ne0$, the transition is a fold: $(F_{\rm M})_{yy}=-f''(a_*)\ne0$, and the branch marker has derivative $f''(a_*)$. Thus (E) holds as $a_0$ varies, while the fixed-point equation remains regular. These are computable recognition conditions for an evolution fold in the machining response. On the regular branch, the multiplier is $f'(a_0)/(1+f'(a_0))$; it reaches $-1$ at slope $-1/2$, the flip threshold studied in \cite{KissBachrathy2022}.

An exact local illustration is the quadratic response
\[
 f(a)=f_*-(a-a_*)+\tfrac12c(a-a_*)^2,
 \qquad f_*>0,\quad c\ne0,
\]
with $a_*$ in the interior of the cutting regime. For
$\lambda=a_0-a_*$, the marker is $c\lambda$ and the forward
multiplier is $1-(c\lambda)^{-1}$. The fixed point
$x=y=f(a_0)$ continues through this evolution fold.

\paragraph{Implicit integration.}
For a scalar differential equation $\dot u=f(u)$ and step size $\delta>0$, backward Euler gives
\begin{equation}\label{eq:application-euler}
 F_\delta(x,y)=y-x-\delta f(y)=0,\qquad A=-1,\quad B=1-\delta f'(y).
\end{equation}
Multiple solutions and turning points of such step equations are studied in \cite{GreenPatrickSpiteri2019}. Let $u_*$ satisfy $f(u_*)=0$, $f'(u_*)>0$ and $f''(u_*)\ne0$. The equilibrium persists for every $\delta$. At $\delta_*=1/f'(u_*)$, one has $B=0$, $(F_\delta)_{yy}=-\delta_*f''(u_*)\ne0$ and $L=-1$. The marker along the equilibrium has derivative $-f'(u_*)$, so (E) holds with parameter $\lambda=\delta-\delta_*$. Its forward multiplier is exactly $-[f'(u_*)\lambda]^{-1}$. This is a singularity of the numerical evolution at a persistent equilibrium of the differential equation.

\paragraph{A mixed interaction in backward Euler.}
For an exact planar example, consider
\begin{equation}\label{eq:euler-mixed-ode}
 \dot u=u-u^3,\qquad \dot v=(\alpha-u)v-v^2,
\end{equation}
with $\alpha$ near $1/2$. Its origin is a hyperbolic equilibrium with eigenvalues $1,\alpha$. A backward-Euler step of size $\delta>0$ has the relation $F_\delta(x,y;\alpha)=H_{\delta,\alpha}(y)-x$, where $x,y\in\R^2$ and
\begin{equation}\label{eq:euler-mixed-backward}
 H_{\delta,\alpha}(u,v)=\bigl((1-\delta)u+\delta u^3,
                      (1-\delta\alpha+\delta u)v+\delta v^2\bigr).
\end{equation}
The origin remains the unique nearby fixed point. At $(\delta,\alpha)=(2,1/2)$,
\[
 A=-I_2,\qquad B=\diag(-1,0),\qquad
 L_+=\diag(-2,-1),\qquad L_-=\diag(0,-1).
\]
With $v_-=w_-=(1,0)^T$, the doubled range reconstruction is $R=(s,0)$, $S=(-s,0)$. It gives
\begin{align*}
 g(s,\delta,\alpha)&=s(\delta-2-\delta s^2),\\
 \sigma_+(s,\delta,\alpha)&=(1-\delta+3\delta s^2)(1-\delta\alpha-\delta s).
\end{align*}
Consequently, at the collision, the coefficients of Section~\ref{sec:mixed-flip} satisfy $a_3=-2$, $c_{\mathrm{PD}}=2$ and
$\det D_{(\delta,\alpha)}\bigl(\delta-2,(1-\delta)(1-\delta\alpha)\bigr)=2$.
The evolution fold coefficient is $4$, so \textup{(PD)} and \textup{(SP)} hold.

The local parameter coordinates
\[
 \mu=\frac{\delta-2}{\delta},\qquad \nu=\frac1\delta-\alpha
\]
put the primary flip and singularity curves at $\mu=0$ and $\nu=0$. For $\mu>0$, the small two-cycle consists of $(\sqrt\mu,0)$ and $(-\sqrt\mu,0)$. Its singularity curve and regular forward multipliers are
\begin{equation}\label{eq:euler-mixed-spectrum}
 \mu=\nu^2,\qquad
 m_{\rm c}=(2\delta-5)^{-2},\qquad
 \rho_\infty=\bigl[\delta^2(\nu^2-\mu)\bigr]^{-1}.
\end{equation}
Indeed, the backward derivative at either state is diagonal, with entries $2\delta-5$ and $\delta(\nu\pm\sqrt\mu)$; multiplication over the two steps gives these formulas. The mixed interaction occurs in the numerical evolution while the underlying differential equation retains its hyperbolic equilibrium.

\section{Further directions}\label{sec:scope}
The joint reduction suggests a programme for higher codimension.
At an interaction with a cusp or a degenerate period doubling,
higher-order coefficients of the reduced periodic equation must be
combined with the transition markers. With several singular steps,
the separate markers record which transitions fail and how their
singularity sets meet. The next task is to derive simultaneous normal
forms and recognition conditions for these configurations, extending
the two-parameter diagrams obtained here.

Two further problems involve different linear degeneracies:
an evolution fold meeting a non-real unit-modulus pair, and
simultaneous singularity of both state derivatives at one transition.
The former calls for a reduction respecting the oscillatory modes;
the latter requires a joint analysis of the forward and backward
transition equations. In each case, the periodic equations and
transition markers specify the geometry to be classified.

\section*{Declarations}
\noindent\textbf{Funding.}
Jorge Buescu was partially supported by Funda\c{c}\~ao para a Ci\^encia e a Tecnologia (FCT), Portugal, through CEMS.UL under project UID/04561/2025 (\href{https://doi.org/10.54499/UID/04561/2025}{doi:10.54499/UID/04561/2025}). Henrique M. Oliveira was partially supported by FCT through CAMGSD---Centro de An\'alise Matem\'atica, Geometria e Sistemas Din\^amicos---under projects UID/04459/2025 (\href{https://doi.org/10.54499/UID/04459/2025}{doi:10.54499/UID/04459/2025}) and UID/PRR/04459/2025 (\href{https://doi.org/10.54499/UID/PRR/04459/2025}{doi:10.54499/UID/PRR/04459/2025}), financed by national funds and by the European Union through the Recovery and Resilience Plan (PRR)---NextGenerationEU.

\medskip
\noindent\textbf{Author contributions.} Both authors contributed equally to all stages of the work.

\medskip
\noindent\textbf{Competing interests.} The authors declare no competing financial or non-financial interests.

\medskip
\noindent\textbf{Data availability.} No empirical datasets were generated or analysed. The proofs and exact calculations are presented in the article.

\medskip
\noindent\textbf{Use of generative AI.} ChatGPT was used for  language and editorial revision. The authors take full responsibility for the mathematical arguments, calculations and final text.

\bibliographystyle{abbrv}
\bibliography{singular_evolution_references}
\end{document}